\documentclass[12pt]{article}
\usepackage{amssymb,amsmath,latexsym, epsfig}

\begin{document}

\begin{doublespace}

\newtheorem{thm}{Theorem}[section]
\newtheorem{lemma}[thm]{Lemma}
\newtheorem{defn}{Definition}[section]
\newtheorem{prop}[thm]{Proposition}
\newtheorem{corollary}[thm]{Corollary}
\newtheorem{remark}[thm]{Remark}
\newtheorem{example}[thm]{Example}
\numberwithin{equation}{section}

\def\ee{\varepsilon}
\def\qed{{\hfill $\Box$ \bigskip}}
\def\MM{{\cal M}}
\def\BB{{\cal B}}
\def\LL{{\cal L}}
\def\FF{{\cal F}}
\def\GG{{\cal G}}
\def\EE{{\cal E}}
\def\QQ{{\cal Q}}

\def\R{{\mathbb R}}
\def\L{{\bf L}}
\def\E{{\mathbb E}}
\def\F{{\bf F}}
\def\P{{\mathbb P}}
\def\N{{\mathbb N}}
\def\eps{\varepsilon}
\def\wh{\widehat}
\def\pf{\noindent{\bf Proof.} }

\title{\Large \bf On the Asymptotic Normality of the Conditional Maximum Likelihood Estimators for the Truncated Regression Model and the Tobit Model         }
\author{ Chunlin Wang\\
Department of Economics\\
University of Pennsylvania\\
Philadelphia, PA 19104\\
 Email: chunlinw@sas.upenn.edu\\   }
\date{}
\maketitle

\begin{abstract}
In this paper, we study the asymptotic normality of the
conditional maximum likelihood (ML) estimators for the truncated
regression model and the Tobit model. We show that under the
general setting assumed in his book, the conjectures made by
Hayashi (2000) \footnote{see page 516, and page 520 of Hayashi
(2000).} about the asymptotic normality of the conditional ML
estimators for both models are true, namely, a sufficient
condition is the nonsingularity of  $\mathbf{x_tx'_t}$.
\end{abstract}



 \vspace{.1truein}
\noindent{\bf Running Title:} Asymptotic normality of the
conditional ML estimators for the truncated regression model and
the Tobit model

\vfill

\pagebreak
\section{Introduction}
The truncated regression model and the Tobit model (censored
regression model) are two important basic models appearing in many
applications in economics and other areas. The method of
conditional maximum likelihood (ML) can be used to estimate the
parameters in both models. In order to apply this method, the
consistency and asymptotic normality of the estimator have to be
verified. In the following, these two models and their conditional
ML estimators are introduced in the exactly same way as what
Hayashi (2000) did.

\subsection{Truncated Regression Model}

 For the truncated regression model, the
following assumptions are adopted:

{\bf Assumption 1}. Suppose that $\{y_t,\mathbf{x}_t\}$ is i.i.d
satisfying
\begin{align}
y_t&=\mathbf{x}'_t\boldsymbol{\beta}_0+\epsilon_t, \\
 \epsilon_t|\mathbf{x}_t&\sim N(0,\sigma^2_{0}), \quad
t=1,2,...,n,
\end{align}
where $\mathbf{x}_t$ and $\boldsymbol{\beta}_0$ are both vectors
with $K$ components.

{\bf Assumption 2}. The truncation rule is: $y_t>c$ where $c$ is a
known constant. Only those observations satisfying the truncation
rule are included in the sample.

Since $y_t|\mathbf{x}_t\sim N(\mathbf{x'_t\boldsymbol{\beta_0}},
\sigma^2_0)$, it can be established that
\begin{align}
E(y_t|\mathbf{x}_t,
y_t>c)&=\mathbf{x}'_t\boldsymbol{\beta_0}+\sigma_0\lambda({ \tfrac{c-\mathbf{x}'_t\boldsymbol{\beta_0}}{\sigma_0}}),\\
Var(y_t|\mathbf{x}_t,
y_t>c)&=\sigma^2_0\left\{1-\lambda(\tfrac{c-\mathbf{x}'_t\boldsymbol{\beta_0}}{\sigma_0})\left[
\lambda(\tfrac{c-\mathbf{x}'_t\boldsymbol{\beta_0}}{\sigma_0})-\tfrac{c-\mathbf{x}'_t\boldsymbol{\beta_0}}{\sigma_0}\right]\right\},
\end{align}
where $\displaystyle
\lambda(\tfrac{c-\mathbf{x}'_t\boldsymbol{\beta_0}}{\sigma_0})\equiv\tfrac{\phi(\tfrac{c-\mathbf{x}'_t\boldsymbol{\beta_0}}{\sigma_0})}
{1-\Phi(\tfrac{c-\mathbf{x}'_t\boldsymbol{\beta_0}}{\sigma_0})}$
with $\phi$ as the density of $N(0,1)$ and $\Phi$ as the
cumulative distribution function of $N(0,1)$. $\lambda$ is also
called the inverse Mill's ratio.

 The log conditional likelihood for observation $t$ is:
\begin{equation}
\log
f(y_t|\mathbf{x}_t;\boldsymbol{\beta},\sigma^2)=\left\{-\frac{1}{2}\log(2\pi)-\frac{1}{2}\log(\sigma^2)-\frac{1}{2}
{\left(\frac
{y_t-\mathbf{x}'_t\boldsymbol{\beta}}{\sigma}\right)}^2\right\}-\log\left[1-\Phi\left(\frac{c-\mathbf{x}'_t\boldsymbol{\beta}}{\sigma}\right)\right],
\end{equation}
where $(\boldsymbol{\beta},\sigma^2)$ are the hypothetical values
of $(\boldsymbol{\beta}_0,\sigma_0^2)$ , and $\Phi$ is the
cumulative distribution function of $N(0,1)$.

For simplification, the following reparameterization is used:
\begin{equation}
\boldsymbol{\delta}=\boldsymbol{\beta}/{\sigma},\quad
\gamma={1}/{\sigma}.
\end{equation}

The reparameterized log conditional likelihood is
\begin{equation}
\log
\tilde{f}(y_t|\mathbf{x}_t;\boldsymbol{\delta},\gamma)=\Bigl[-\frac{1}{2}\log(2\pi)+\log(\gamma)-\frac{1}{2}
{(\gamma
y_t-\mathbf{x}'_t\boldsymbol{\delta})}^2\Bigr]-\log[1-\Phi(\gamma
c-\mathbf{x}'_t\boldsymbol{\delta})].
\end{equation}

The objective function in ML estimation is the average log
conditional likelihood of the sample. The conditional ML estimator
$(\hat{\boldsymbol{\delta}},\hat{\gamma})$ of
$({\boldsymbol{\delta_0}},\hat{\gamma}_0)$ is the
$({\boldsymbol{\delta}},{\gamma})$ that maximizes the objective
function.

Hayashi (2000) gave the following expressions of the score and the
Hessian for observation $t$:
\begin{align}
\mathbf{s}(\mathbf{w}_t;\boldsymbol{\delta},\gamma)&=\begin{bmatrix}
(\gamma
y_t-\mathbf{x}'_t\boldsymbol{\delta})\mathbf{x}_t\\
\frac{1}{\gamma}-(\gamma y_t-\mathbf{x}'_t\boldsymbol{\delta})y_t
\end{bmatrix}+\lambda(v_t)\begin{bmatrix}-\mathbf{x}_t\\
c\end{bmatrix},
\\
\notag\\
\bigskip
\mathbf{H}(\mathbf{w}_t;\boldsymbol{\delta},\gamma)&=-\begin{bmatrix}
\mathbf{x}_t\mathbf{x}'_t&-y_t\mathbf{x}_t\\
-y_t\mathbf{x}'_t&\frac{1}{\gamma^2}+y^2_t
\end{bmatrix}+\lambda(v_t)[\lambda(v_t)-v_t]\begin{bmatrix}
\mathbf{x}_t\mathbf{x}'_t&-c\mathbf{x}_t\\
-c\mathbf{x}'_t&c^2
\end{bmatrix},
\end{align}\\
where $\mathbf{s}(\mathbf{w}_t;\boldsymbol{\delta},\gamma)$ is a
vector of dimension $(K+1)\times 1$,
$\mathbf{H}(\mathbf{w}_t;\boldsymbol{\delta},\gamma)$ is a square
matrix of dimension $(K+1)\times (K+1)$ with $K$ as the number of
regressors, $\mathbf{w}_t=(y_t,\mathbf{x}'_t)'$,\,
$\lambda(v_t)\equiv \tfrac{\phi(v_t)}{1-\Phi(v_t)}$ with
$v_t\equiv \gamma
c-\mathbf{x}'_t\boldsymbol{\delta}=\frac{c-\mathbf{x}'_t\boldsymbol{\beta}}{\sigma}$.

By verifying that the conditions of Proposition 1.1 (see below)
are satisfied, Hayashi (2000) proved that the ML estimator
$(\hat{\boldsymbol{\delta}},\hat{\gamma})$ of
$({\boldsymbol{\delta_0}},{\gamma_0})$  is consistent under the
nonsingularity of  $\mathbf{
E(}\mathbf{x}_t\mathbf{x}'_t\mathbf{)}$.

As to asymptotic normality, Hayashi (2000) pointed out that
$E[\mathbf{s}|\mathbf{x}_t]=\mathbf{0}$ and the conditional
information equality holds, i.e.
$E[\mathbf{ss'}|\mathbf{x}_t]=-E[\mathbf{H}|\mathbf{x}_t]$. So
condition $3$ of Proposition 1.3 is satisfied. However conditions
$4$ and $5$ of proposition 1.3  are not verified. For the case
where $\mathbf{\{x_t\}}$ is a sequence of fixed constants, Sapra
(1992) showed asymptotic normality under the assumption that
$\mathbf{x}_t$ is bounded, $\lim_{n \rightarrow
\infty}\frac{1}{n}\sum_{t=1}^{n}\mathbf{x}_t\mathbf{x}'_t$ is
nonsingular and observations are serially correlated.

Hayashi (2000) conjectured that for the case where $\mathbf{x}_t$
is random as in the current setting (i.e., Assumption 1 and 2 are
satisfied), a sufficient condition for asymptotic normality is the
nonsingularity of $\mathbf{
E(}\mathbf{x}_t\mathbf{x}'_t\mathbf{)}$\footnote{see page 516 of
Hayashi (2000).}.

\subsection{Tobit Model}
For the Tobit model, the following assumption is adopted:

{\bf Assumption 1'}. Suppose that $\{y_t,\mathbf{x}_t\}$ is i.i.d
satisfying
\begin{align}
y^*_t&=\mathbf{x}'_t\boldsymbol{\beta_0}+\epsilon_t, \\
\epsilon_t|\mathbf{x}_t&\sim N(0,\sigma^2_{0}),\quad
t=1,2,...,n,\\
 y_t&=\begin{cases}
 \,\,\,\,\textrm{  }y^*_t &\textrm{ }\textrm{if } y^*_t>c,\\
 \,\,\,\,\textrm{  }c    &\textrm{ }\textrm{if } y^*_t\leq c,
 \end{cases}
 \end{align}
where $\mathbf{x}_t$ and $\boldsymbol{\beta}_0$ are both vectors
with $K$ components, $c$ is a known constant. Different from the
truncated regression model in above, here the observations for
which the value of the dependent variable $y^*_t$ doesn't meet the
rule $y^*_t>c$ are included in the sample. Another way to write
the Tobit model is
\begin{equation}
y_t=\max\{\mathbf{x}'_t\boldsymbol{\beta}_0+\epsilon_t,c\}.
\end{equation}

The log conditional likelihood for observation $t$ is:
\begin{equation}
\log
f(y_t|\mathbf{x}_t;\boldsymbol{\beta},\sigma^2)=(1-D_t)\log\left[\frac{1}{\sigma}\phi\left(\frac{y_t-\mathbf{x}'_t\boldsymbol{\beta}}{\sigma}\right)\right]
+D_t\log\Phi\left(\frac{c-\mathbf{x}'_t\boldsymbol{\beta}}{\sigma}\right),
\end{equation}
where $(\boldsymbol{\beta},\sigma^2)$ are the hypothetical values
of $(\boldsymbol{\beta}_0,\sigma_0^2)$, $\phi$ is the density of
$N(0,1)$ and $\Phi$ is the cumulative distribution of $N(0,1)$,
and the dummy variable $D_t$ is defined as
\begin{equation}
D_t=\begin{cases}
    \,\,\,\,\textrm{ }0 &\textrm{ } \textrm{if }y_t>c \,\,(\textrm{i.e.},\, y^*_t>c)\\
    \,\,\,\,\textrm{ }1 &\textrm{ } \textrm{if }y_t=c \,\,(\textrm{i.e.}, \,y^*_t\leq c).
    \end{cases}
\end{equation}\\
The objective function in ML estimation is the average log
conditional likelihood of the sample.

As in the truncation regression model in above, to make analysis
easier, the reparameterization (1.6) is used and the
reparameterized log conditional likelihood is:
\begin{equation}
\log
\tilde{f}(y_t|\mathbf{x}_t;\boldsymbol{\delta},\gamma)=(1-D_t)\Bigl\{-\frac{1}{2}\log(2\pi)+\log(\gamma)-\frac{1}{2}
{(\gamma
y_t-\mathbf{x}'_t\boldsymbol{\delta})}^2\Bigr\}-D_t\log\Phi(\gamma
c-\mathbf{x}'_t\boldsymbol{\delta}).
\end{equation}

Hayashi (2000) gave the following expressions of the score and the
Hessian for observation $t$:
\begin{align}
\displaystyle
\mathbf{s}(\mathbf{w}_t;\boldsymbol{\delta},\gamma)&=(1-D_t)\begin{bmatrix}
(\gamma
y_t-\mathbf{x}'_t\boldsymbol{\delta})\mathbf{x}_t\\
\frac{1}{\gamma}-(\gamma y_t-\mathbf{x}'_t\boldsymbol{\delta})y_t
\end{bmatrix}+D_t\lambda(-v_t)\begin{bmatrix}-\mathbf{x}_t\\
c\end{bmatrix},
\\
\notag\\
\medskip
\mathbf{H}(\mathbf{w}_t;\boldsymbol{\delta},\gamma)&=-(1-D_t)\begin{bmatrix}
\mathbf{x}_t\mathbf{x}'_t&-y_t\mathbf{x}_t\\
-y_t\mathbf{x}'_t&\frac{1}{\gamma^2}+y^2_t
\end{bmatrix}-D_t\lambda(-v_t)[\lambda(-v_t)+v_t]\begin{bmatrix}
\mathbf{x}_t\mathbf{x}'_t&-c\mathbf{x}_t\\
-c\mathbf{x}'_t&c^2
\end{bmatrix},
\end{align}\\
where $\mathbf{s}(\mathbf{w}_t;\boldsymbol{\delta},\gamma)$ is a
vector of dimension $(K+1)\times 1$,
$\mathbf{H}(\mathbf{w}_t;\boldsymbol{\delta},\gamma)$ is a square
matrix of dimension $(K+1)\times (K+1)$, $K$ is the number of
regressors, $\mathbf{w}_t=(y_t,\mathbf{x}'_t)'$,\,
$\lambda(-v_t)\equiv \tfrac{\phi(-v_t)}{1-\Phi(-v_t)}$ with
$v_t\equiv \gamma
c-\mathbf{x}'_t\boldsymbol{\delta}=\frac{c-\mathbf{x}'_t\boldsymbol{\beta}}{\sigma}$.

For consistency, Hayashi (2000) pointed out that the relevant
consistency theorem for the Tobit Model is Proposition 1.2 (see
below)\footnote{see page 520 of Hayashi (2000)}. He also mentioned
that when $\{y_t,\mathbf{x}_t\}$ is ergodic stationary but not
necessary i.i.d, the conditional ML estimator
$(\hat{\boldsymbol{\delta}},\hat{\gamma})$ of
$({\boldsymbol{\delta_0}},{\gamma_0})$ is consistent\footnote{see
exercise 3 on page 521 of Hayashi (2000)}.

For the case where $\mathbf{\{x_t\}}$ is a sequence of fixed
constants, Amemiya (1973) proved the consistency and asymptotic
normality of the conditional ML estimator for the Tobit model
under the assumption that $\mathbf{x}_t$ is bounded and $\lim_{n
\rightarrow
\infty}\frac{1}{n}\sum_{t=1}^{n}\mathbf{x}_t\mathbf{x}'_t$ is
nonsingular.

 For the Tobit model, Hayashi
(2000) conjectured that for the case where $\mathbf{x}_t$ is
random as in the current setting (i.e., Assumption 1' is
satisfied), a sufficient condition for asymptotic normality is the
nonsingularity of $\mathbf{
E(}\mathbf{x}_t\mathbf{x}'_t\mathbf{)}$\footnote{see page 520 of
Hayashi (2000).}.

In this paper, we show that Hayashi's conjectures for the
asymptotic normality of both models are true, i.e. a sufficient
condition for the asymptotic normality for both models is the
nonsingularity of
 $\mathbf{
E(}\mathbf{x}_t\mathbf{x}'_t\mathbf{)}$.

The content of this paper is organized as follows. First in below
we cite three propositions from Hayahsi (2000), which will be used
to show the main results. Then in Section 2, we show that a
sufficient condition for the asymptotic normality for the
truncated regression model is the nonsingularity of $\mathbf{
E(}\mathbf{x}_t\mathbf{x}'_t\mathbf{)}$. In Section 3, we show
that a sufficient condition for the asymptotic normality for the
Tobit model is the nonsingularity of $\mathbf{
E(}\mathbf{x}_t\mathbf{x}'_t\mathbf{)}$

Now we present three propositions from Hayashi (2000).
\begin{prop}
\footnote{see Proposition 7.5 on page 464 of Hayashi (2000). }
{\bf (Consistency of conditional ML with compact parameter
space):} Let $\{y_t,\mathbf{x}_t\}$ be ergodic stationary with
conditional density $f(y_t|\mathbf{x}_t;\boldsymbol{\theta}_0)$
and let $\boldsymbol{\hat{\theta}}$ be the conditional $ML$
estimator, which maximizes the average log conditional likelihood
(derived under the assumption that $\{y_t,\mathbf{x}_t\}$ is
i.i.d.):
$$\displaystyle
\boldsymbol{\hat{\theta}}={\verb"argmax"}_{\boldsymbol{\theta}\in
\Theta}\,\frac{1}{n}\sum^{n}_{t=1}\log
f(y_t|\mathbf{x}_t;\boldsymbol{\theta}).$$ Suppose the model is
correctly specified so that $\boldsymbol{\theta}_0$ is in
$\Theta$. Suppose that (i) the parameter space $\Theta$ is a
compact subset of $\R^p$, (ii) $f(y_t|\mathbf{x}_t;
\boldsymbol{\theta})$ is continuous in $\boldsymbol{\theta}$ for
all $(y_t,\mathbf{x}_t)$, and (iii) $f(y_t|\mathbf{x}_t;
\boldsymbol{\theta})$ is measurable in $(y_t, \mathbf{x}_t)$ for
all $\theta \in \Theta$ (so $\boldsymbol{\hat{\theta}}$ is a
well-defined random variable). Suppose, further that
\begin{enumerate}
\item (identification) $Prob[f(y_t|\mathbf{x}_t;
\boldsymbol{\theta})\neq
f(y_t|\mathbf{x}_t;\boldsymbol{\theta}_0)]>0$ for all
$\boldsymbol{\theta}\neq \boldsymbol{\theta}_0$ in $\Theta$,

\item(dominance) $E[\sup_{\boldsymbol{\theta}\in \Theta}|\log
f(y_t|\mathbf{x}_t;\boldsymbol{\theta})|]<\infty$ (note: the
expectation is over $y_t$ and $\mathbf{x}_t$).
\end{enumerate}
Then $\boldsymbol{\hat{\theta}}\,{\rightarrow}_p\,\,
\boldsymbol{\theta}_0$.
\end{prop}

\begin{prop}
\footnote{see Proposition 7.6 on page 464-465 of Hayashi (2000). }
{\bf (Consistency of conditional ML without compactness):} Let
$\{y_t,\mathbf{x}_t\}$ be ergodic stationary with conditional
density $f(y_t|\mathbf{x}_t;\boldsymbol{\hat{\theta}}_0)$ and let
$\hat{\theta}$ be the conditional $ML$ estimator, which maximizes
the average log conditional likelihood (derived under the
assumption that $\{y_t,\mathbf{x}_t\}$ is i.i.d.):
$$\displaystyle
\boldsymbol{\hat{\theta}}={\verb"argmax"}_{\boldsymbol{\theta}\in
\Theta}\,\frac{1}{n}\sum^{n}_{t=1}\log
f(y_t|\mathbf{x}_t;\boldsymbol{\theta}).$$ Suppose the model is
correctly specified so that $\boldsymbol{\theta}_0$ is in
$\Theta$. Suppose that (i) the true parameter vector
$\boldsymbol{\theta}_0$ is an element of the interior of a convex
parameter space $\Theta$ ($\subset \R^p$), (ii) $\log
f(y_t|\mathbf{x}_t; \boldsymbol{\theta})$ is concave in
$\boldsymbol{\theta}$ for all $(y_t,\mathbf{x}_t)$, and (iii)
$\log f(y_t|\mathbf{x}_t; \boldsymbol{\theta})$ is measurable in
$(y_t, \mathbf{x}_t)$ for all $\theta \in \Theta$. (For
sufficiently large $n$, $\boldsymbol{\hat{\theta}}$ well-defined).
Suppose, further that
\begin{enumerate}
\item (identification) $Prob[f(y_t|\mathbf{x}_t;
\boldsymbol{\theta})\neq
f(y_t|\mathbf{x}_t;\boldsymbol{\theta}_0)]>0$ for all
$\boldsymbol{\theta}\neq \boldsymbol{\theta}_0$ in $\Theta$,

\item $E[|\log f(y_t|\mathbf{x}_t;\boldsymbol{\theta})|]<\infty$
(i.e., $E[\log f(y_t|\mathbf{x}_t;\boldsymbol{\theta})]$ exists
and is finite) for all $\boldsymbol{\theta}\in {\Theta}$(note: the
expectation is over $y_t$ and $\mathbf{x}_t$).
\end{enumerate}
Then as $n \rightarrow \infty$, $\boldsymbol{\hat{\theta}}$ exists
with probability approaching $1$ and
$\boldsymbol{\hat{\theta}}\,{\rightarrow}_p\,\,
\boldsymbol{\theta}_0$.
\end{prop}

\begin{prop}\footnote{see Proposition 7.9 on page 475 of Hayashi
(2000).} {\bf (Asymptotic normality of conditional ML):} Let
$\mathbf{w}_t$ ($\equiv (y_t,\mathbf{x}'_t)'$) be i.i.d. Suppose
the conditions of either Proposition 1.1 or Proposition 1.2 are
satisfied, so that $\boldsymbol{\hat{\theta}}\,{\rightarrow}_p\,\,
\boldsymbol{\theta}_0$. Suppose, in addition, that
\begin{enumerate}
\item $\boldsymbol{\theta}_0$ is in the interior of $\Theta$,
\item $f(y_t|\mathbf{x}_t;\boldsymbol{\theta})$ is twice
continuously differentiable in $\boldsymbol{\theta}$ for all
$(y_t,\mathbf{x}_t)$,

\item
$E[\mathbf{s}(\mathbf{w}_t;\boldsymbol{\theta}_0)]=\mathbf{0}$ and
$-E[\mathbf{H}(\mathbf{w}_t;\boldsymbol{\theta}_0)]=E[\mathbf{s}(\mathbf{w}_t;\boldsymbol{\theta}_0)
\mathbf{s}(\mathbf{w}_t;\boldsymbol{\theta}_0)']$, where
$\mathbf{s}$ and $\mathbf{H}$ functions are the score and the
Hessian for observation $t$.

\item (local dominance condition on the Hessian) for some
neighborhood ${\cal N}$ of $\boldsymbol{\theta}_0$,
$$ E[\sup_{\boldsymbol{\theta}\in {\cal
N}}\|\mathbf{H}(\mathbf{w}_t;\boldsymbol{\theta})\|]<\infty,$$ so
that for any consistent estimator $\boldsymbol{\tilde{\theta}}$,
$\frac{1}{n}\sum^{n}_{t=1}\mathbf{H}(\mathbf{w}_t;\boldsymbol{\tilde{\theta}})\,{\rightarrow}_p\,
E[\mathbf{H}(\mathbf{w}_t;\boldsymbol{\theta}_0)],$ \item
$E[\mathbf{H}(\mathbf{w}_t;\boldsymbol{\theta}_0)]$ is
nonsingular.
\end{enumerate}
Then $\boldsymbol{\hat{\theta}}$ is asymptotic normal with
$\verb"Avar"(\boldsymbol{\hat{\theta}})$ given by the following:
$$\verb"Avar"(\boldsymbol{\hat{\theta}})=-\{E[\mathbf{H}(\mathbf{w}_t;\boldsymbol{\theta}_0)]\}^{-1}=\{E[\mathbf{s}(\mathbf{w}_t;\boldsymbol{\theta}_0)
\mathbf{s}(\mathbf{w}_t;\boldsymbol{\theta}_0)']\}^{-1}.$$
\end{prop}

\section{Asymptotic Normality of the Conditional ML of the Truncated Regression Model}
As we mentioned in above, for the truncated regression model
introduced in subsection 1.1, Hayashi (2000) showed that the
conditional ML estimator satisfies the conditions of Proposition
1.1 under the nonsingularity of $\mathbf{
E(}\mathbf{x}_t\mathbf{x}'_t\mathbf{)}$. Therefore, by Proposition
1.1, the conditional ML estimator
$(\boldsymbol{\hat{\delta}},\hat{\gamma})$ of
$(\boldsymbol{\delta}_0,{\gamma}_0)$ is consistent when $\mathbf{
E(}\mathbf{x}_t\mathbf{x}'_t\mathbf{)}$ is nonsingular. For
asymptotic normality, he mentioned that condition $3$ of
Proposition 1.3 is satisfied. It is easy to see that conditions
$1$ and $2$ of Proposition 1.3 are satisfied.

In this section, we show the following theorem holds.

\begin{thm}
 For the truncated regression model satisfying Assumptions 1 and
 2, if $\mathbf{
E(}\mathbf{x}_t\mathbf{x}'_t\mathbf{)}$ is nonsingular, then
conditions $4$ and $5$ of Proposition 1.3 are satisfied.
\end{thm}

\pf First we show that condition $4$ of Proposition 1.3 is
satisfied.

 Define $\mathbf{A}$ as the matrix
\begin{align}
-\begin{bmatrix}
\mathbf{x_tx'_t}&-y_t\mathbf{x}_t\\
-y_t\mathbf{x}'_t&\frac{1}{\gamma^2}+y^2_t
\end{bmatrix},
\end{align}
and define $\mathbf{B}$ as the matrix
\begin{align}
\lambda(v_t)[\lambda(v_t)-v_t]\begin{bmatrix}
\mathbf{x}_t\mathbf{x}'_t&-c\mathbf{x}_t\\
-c\mathbf{x}'_t&c^2
\end{bmatrix}.
\end{align}
By the expression of the Hessian
$\mathbf{H}(\mathbf{w}_t;\boldsymbol{\delta},\gamma)$ in (1.9),
\begin{equation}
\mathbf{H}(\mathbf{w}_t;\boldsymbol{\delta},\gamma)=\mathbf{A}+\mathbf{B}.
\end{equation}
Therefore,
\begin{equation}
\|\mathbf{H}(\mathbf{w}_t;\boldsymbol{\delta},\gamma)\|\leq
\|\mathbf{A}\|+\|\mathbf{B}\|,
\end{equation}
where $\|\cdot\|$ is the Euclidean norm of a matrix, which is
defined as the square root of the sum of squares of the elements
of the matrix.

It is easy to see that
\begin{equation}
{\|\mathbf{A}\|}^2={\|\mathbf{x_tx'_t}\|}^2+2{\|y_t\mathbf{x}_t\|}^2+{\displaystyle
(1/{\gamma^2}+y^2_t)}^2.
\end{equation}
Since $\mathbf{x}_t$ is a vector of $K$ components, we write it as
$(\mathbf{x}_{t1},\mathbf{x}_{t2},...,\mathbf{x}_{tK})'.$ We have
\begin{equation}
\displaystyle {\|y_t\mathbf{x}_t\|}^2\leq
\sum_{i=1}^{K}\tfrac{1}{2}({y^4_t}+{\mathbf{x}^4_{ti}})\leq
\tfrac{K}{2}y^4_t+\tfrac{1}{2}{\|\mathbf{x}_t\mathbf{x}'_t\|}^2.
\end{equation}
Thus
\begin{align}
{\|\mathbf{A}\|}^2&\leq
{\|\mathbf{x}_t\mathbf{x}'_t\|}^2+2\tfrac{K}{2}y^4_t+{\|\mathbf{x}_t\mathbf{x}'_t\|}^2+2/\gamma^4+2y^4_t \notag\\
&\leq2{\|\mathbf{x}_t\mathbf{x}'_t\|}^2+(K+2)y^4_t+2/\gamma^4,
\end{align}
which implies
\begin{equation}
{\|\mathbf{A}\|}\leq
\sqrt{2}\,{\|\mathbf{x}_t\mathbf{x}'_t\|}+\sqrt{K+2}\,y^2_t+\sqrt{2}/\gamma^2.
\end{equation}
By (1.1) and (1.6), we have
\begin{equation}
y^2_t\leq
\tfrac{2}{\gamma_0^2}{(\mathbf{x}'_t{\boldsymbol{\delta}_0})}^2+2{\epsilon^2_t}.
\end{equation}

 Since $\boldsymbol{\delta}$ and $\boldsymbol{\delta}_0$ are
vectors of $K$ components, we write them as
$(\boldsymbol{\delta}_1,...,\boldsymbol{\delta}_K)'$ and
$(\boldsymbol{\delta}_{01},...,$\\$\boldsymbol{\delta}_{0K})'$.
Define the neighborhood $\cal N$ of $(\boldsymbol{\delta}_0,
{\gamma}_0)$ ($\equiv \boldsymbol{\theta}_0$) as
$$\{(\boldsymbol{\delta}, {\gamma}):
\max_{i=1,...,K}|\boldsymbol{\delta}_i-\boldsymbol{\delta}_{0i}|<
C_1,\, |\gamma-\gamma_0|< C_2,\textrm{ with } C_2 \textrm{
satisfies } 0<-C_2+\gamma_0 \},$$

where $C_1$ and $C_2$ are positive constants.

 Since
$\gamma_0=1/\sigma_0$ and $\sigma_0$ is finite, $\gamma_0\neq 0$.
So we can always find a small positive $C_2$ such that
$0<-C_2+\gamma_0$. This implies for any
$(\boldsymbol{\delta,\gamma})\in {\cal N}$, $1/\gamma\leq
\frac{1}{-C_2+\gamma_0}$. Combining this with (2.8) and (2.9), we
have, for any $(\boldsymbol{\delta},\gamma)\in {\cal N}$,
\begin{align}
\sup_{(\boldsymbol{\delta},\gamma)\in {\cal N}}\|\mathbf{A}\|&\leq
\sqrt{2}\,{\|\mathbf{x}_t\mathbf{x}'_t\|}+\sqrt{K+2}\,\tilde{C}\|\mathbf{x}_t\mathbf{x}'_t\|+
2\sqrt{K+2}\,{\epsilon^2_t}+\sqrt{2}\,\tfrac{1}{{(-C_2+\gamma_0)}^2}\notag \\
&\leq
(\sqrt{2}+\sqrt{K+2}\,\tilde{C})\|\mathbf{x}_t\mathbf{x}'_t\|+2\sqrt{K+2}\,{\epsilon^2_t}+\sqrt{2}\,\tfrac{1}{{(-C_2+\gamma_0)}^2},
\end{align}
where
$\tilde{C}=\frac{2}{\gamma_0^2}\max(\delta^2_{01},...,\delta^2_{0K})$.

 Next we look at the matrix $\mathbf{B}$ defined in (2.2).

 It is well known that as the derivative of
 $\lambda(v_t)$, $\lambda'(v_t)$ satisfies
 \begin{equation}
 \lambda'(v_t)=\lambda(v_t)(\lambda(v_t)-v_t),
 \end{equation}
 and $\lambda(v_t)$ is between $0$ and $1$. Therefore
 \begin{equation}
 {\|\mathbf{B}\|}^2\leq
 {\|\mathbf{x_tx'_t}\|}^2+2{\|c\mathbf{x}_t\|}^2+c^4.
 \end{equation}
 Since ${\|c\mathbf{x}_t\|}^2$ satisfies
\begin{equation}
\displaystyle {\|c\mathbf{x}_t\|}^2\leq
\sum_{i=1}^{K}\tfrac{1}{2}({c^4}+{\mathbf{x}^4_{ti}})\leq
\tfrac{K}{2}c^4+\tfrac{1}{2}{\|\mathbf{x}_t\mathbf{x}'_t\|}^2,
\end{equation}
\begin{equation}
 {\|\mathbf{B}\|}^2\leq
 {2\|\mathbf{x}_t\mathbf{x}'_t\|}^2+Kc^4+c^4,
 \end{equation}
which implies
\begin{equation}
{\|\mathbf{B}\|}\leq
\sqrt{2}\|\mathbf{x}_t\mathbf{x}'_t\|+\sqrt{K+1}c^2.
\end{equation}
Combining this with (2.10) and (2.4), we have
\begin{equation}
\sup_{(\boldsymbol{\delta},\gamma)\in {\cal
N}}\|\mathbf{H}(\mathbf{w}_t;\boldsymbol{\delta},\gamma)\|\leq
(2\sqrt{2}+\sqrt{K+2}\,\tilde{C})\|\mathbf{x}_t\mathbf{x}'_t\|+
2\sqrt{K+2}\,{\epsilon^2_t}+\sqrt{2}\,\tfrac{1}{{(-C_2+\gamma_0)}^2}+
\sqrt{K+1}c^2.
\end{equation}
Since
$E[\epsilon^2_t|\mathbf{x}_t]=\sigma^2_0=\frac{1}{\gamma^2_0}$,
\begin{align}
E[\sup_{(\boldsymbol{\delta},\gamma)\in {\cal
N}}\|\mathbf{H}(\mathbf{w}_t;\boldsymbol{\delta},\gamma)\|]&\leq
(2\sqrt{2}+\sqrt{K+2}\,\tilde{C})E[\|\mathbf{x}_t\mathbf{x}'_t\|]+
2\sqrt{K+2}\,\frac{1}{\gamma^2_0}+\sqrt{2}\,\tfrac{1}{{(-C_2+\gamma_0)}^2}\notag\\
&\quad+\sqrt{K+1}c^2.
\end{align}
We know that $E[\|\mathbf{x}_t\mathbf{x}'_t\|]<\infty$ if
$E[\mathbf{x}_t\mathbf{x}'_t]$ exists and is finite, and
$E[\mathbf{x}_t\mathbf{x}'_t]$ exists and is finite if
$E[\mathbf{x}_t\mathbf{x}'_t]$ is nonsingular. Therefore, when
$E[\mathbf{x}_t\mathbf{x}'_t]$ is nonsingular,
\begin{equation}
E[\sup_{(\boldsymbol{\delta},\gamma)\in {\cal
N}}\|\mathbf{H}(\mathbf{w}_t;\boldsymbol{\delta},\gamma)\|]<
\infty,
\end{equation}
namely, condition $4$ of Proposition 1.3 is satisfied.\\

Next we show that condition 5 of Proposition 1.3 is satisfied.

Since condition 3 of Proposition 1.3 is satisfied,
\begin{equation}
-E[\mathbf{H}(\mathbf{w}_t;\boldsymbol{\delta}_0,
\gamma_0)]=E[\mathbf{s}(\mathbf{w}_t;\boldsymbol{\delta}_0,
\gamma_0) \mathbf{s}(\mathbf{w}_t;\boldsymbol{\delta}_0,
\gamma_0)'].
 \end{equation}

 It is clear that
$\mathbf{s}(\mathbf{w}_t;\boldsymbol{\delta}_0, \gamma_0)
\mathbf{s}(\mathbf{w}_t;\boldsymbol{\delta}_0, \gamma_0)'$ is
positive semidefinite. This implies that \\$\displaystyle
E[\mathbf{s}(\mathbf{w}_t;\boldsymbol{\delta}_0,
\gamma_0)\mathbf{s}(\mathbf{w}_t;\boldsymbol{\delta}_0,
\gamma_0)']$ is positive semidefinite. Therefore
$E[\mathbf{H}(\mathbf{w}_t;\boldsymbol{\delta}_0, \gamma_0)]$ is
negative semidefinite.

Let
$\mathbf{z}\equiv(\mathbf{z}_1,...,\mathbf{z}_K,\mathbf{z}_{K+1})'\in
\R^{K+1}$ be a solution to the equation
\begin{equation}\label{A}
      {\mathbf{z}}'E[\mathbf{H}(\mathbf{w}_t;\boldsymbol{\delta}_0,
\gamma_0)]\mathbf{z}=0.
\end{equation}
We know that if $\mathbf{z}=(0,...,0)$ is the only solution to the
above equation, then
$E[\mathbf{H}(\mathbf{w}_t;\boldsymbol{\delta}_0, \gamma_0)]$ is
nonsingular.

By (2.19), (2.20) is equivalent to
\begin{equation}
{\mathbf{z}}'E[\mathbf{s}(\mathbf{w}_t;\boldsymbol{\delta}_0,
\gamma_0) \mathbf{s}(\mathbf{w}_t;\boldsymbol{\delta}_0,
\gamma_0)'] \mathbf{z}=0,
\end{equation}
namely,
\begin{equation}
E[{\mathbf{z}}'\mathbf{s}(\mathbf{w}_t;\boldsymbol{\delta}_0,
\gamma_0) \mathbf{s}(\mathbf{w}_t;\boldsymbol{\delta}_0,
\gamma_0)'\mathbf{z}]=0,
\end{equation}
because $\mathbf{z}$ is not random.

In terms of the expression of $\mathbf{s}$ in (1.8), (2.22) is
equivalent to
\begin{equation}
E\biggl[\sum_{i=1}^{K}\bigl(\gamma_0y_t-\mathbf{x}'_t\boldsymbol{\delta}_0-\lambda(v_{0t})\bigr)\mathbf{x}_{ti}\mathbf{z}_i1_{\{y_t>c\}}
+\bigl(\frac{1}{\gamma_0}-(\gamma_0y_t-\mathbf{x}'_t\boldsymbol{\delta}_0)y_t+\lambda(v_{0t})c\bigr)\mathbf{z}_{K+1}1_{\{y_t>c\}}\biggr]^2=0,
\end{equation}
where
$v_{0t}\equiv\gamma_0y_t-\mathbf{x}'_t\boldsymbol{\delta}_0$.

This implies
\begin{equation}
\sum_{i=1}^{K}\bigl(\gamma_0y_t-\mathbf{x}'_t\boldsymbol{\delta}_0-\lambda(v_{0t})\bigr)\mathbf{x}_{ti}\mathbf{z}_i1_{\{y_t>c\}}
+\bigl(\frac{1}{\gamma_0}-(\gamma_0y_t-\mathbf{x}'_t\boldsymbol{\delta}_0)y_t+\lambda(v_{0t})c\bigr)\mathbf{z}_{K+1}1_{\{y_t>c\}}=0
 \,\,\textrm{ a.e.},
\end{equation}
where a.e. means almost everywhere.

The left side of (2.24)
\begin{align}
      =&\sum_{i=1}^{K}({\gamma}_0y_t-\mathbf{x}'_t\boldsymbol{\delta}_0)\mathbf{x}_{ti}\mathbf{z}_i1_{\{y_t>c\}}-\sum_{i=1}^{K}\lambda(v_{0t})\mathbf{x}_{ti}\mathbf{z}_i1_{\{y_t>c\}}\notag\\
       &+\bigl(\frac{1}{{\gamma}_0}-(\gamma_0y_t-\mathbf{x}'_t\boldsymbol{\delta}_0)\frac{1}{\gamma_0}(\gamma_0 y_t-\mathbf{x}'_t\boldsymbol{\delta}_0+
       \mathbf{x}'_t\boldsymbol{\delta}_0)+\lambda(v_{0t})c\bigr)\mathbf{z}_{K+1}1_{\{y_t>c\}}\notag\\
       =&\sum_{i=1}^{K}({\gamma}_0y_t-\mathbf{x}'_t\boldsymbol{\delta}_0)\mathbf{x}_{ti}\mathbf{z}_i1_{\{y_t>c\}}
      -\frac{1}{\gamma_0}(\gamma_0
      y_t-\mathbf{x}'_t\boldsymbol{\delta}_0)\mathbf{x}'_t\boldsymbol{\delta}_0\mathbf{z}_{K+1}1_{\{y_t>c\}}\notag\\
      &-\frac{1}{\gamma_0}(\gamma_0
      y_t-\mathbf{x}'_t\boldsymbol{\delta}_0)^2\mathbf{z}_{K+1}1_{\{y_t>c\}}-\sum_{i=1}^{K}\lambda(v_{0t})\mathbf{x}_{ti}\mathbf{z}_i1_{\{y_t>c\}}
      +\bigl(\frac{1}{\gamma_0}+\lambda(v_{0t})c\bigr)\mathbf{z}_{K+1}1_{\{y_t>c\}}\notag
\end{align}
\begin{align}
      =&({\gamma}_0y_t-\mathbf{x}'_t\boldsymbol{\delta}_0)\biggl[\sum_{i=1}^{K}\mathbf{x}_{ti}\mathbf{z}_i
      -\frac{1}{\gamma_0}\mathbf{x}'_t\boldsymbol{\delta}_0\mathbf{z}_{K+1}\biggr]1_{\{y_t>c\}}\notag\\
      &-\frac{1}{\gamma_0}(\gamma_0
      y_t-\mathbf{x}'_t\boldsymbol{\delta}_0)^2\mathbf{z}_{K+1}1_{\{y_t>c\}}+\biggl[-\sum_{i=1}^{K}\lambda(v_{0t})\mathbf{x}_{ti}\mathbf{z}_i
      +\bigl(\frac{1}{\gamma_0}+\lambda(v_{0t})c\bigr)\mathbf{z}_{K+1}\biggr]1_{\{y_t>c\}}\notag\\
      =&\gamma_0\epsilon_t\biggl[\sum_{i=1}^{K}\mathbf{x}_{ti}\mathbf{z}_i
      -\frac{1}{\gamma_0}\mathbf{x}'_t\boldsymbol{\delta}_0\mathbf{z}_{K+1}\biggr]1_{\{\epsilon_t>c-\frac{1}{\gamma_0}\mathbf{x}'_t\boldsymbol{\delta}_0\}}
      -\gamma_0\epsilon_t^2\mathbf{z}_{K+1}1_{\{\epsilon_t>c-\frac{1}{\gamma_0}\mathbf{x}'_t\boldsymbol{\delta}_0\}}\notag\\
      &+\biggl[-\sum_{i=1}^{K}\lambda(v_{0t})\mathbf{x}_{ti}\mathbf{z}_i
      +\bigl(\frac{1}{\gamma_0}+\lambda(v_{0t})c\bigr)\mathbf{z}_{K+1}\biggr]1_{\{\epsilon_t>c-\frac{1}{\gamma_0}\mathbf{x}'_t\boldsymbol{\delta}_0\}}.
\end{align}

Define $f(\mathbf{x}_t,\mathbf{z},\boldsymbol{\delta}_0,\gamma_0)$
as $\displaystyle \biggl[\sum_{i=1}^{K}\mathbf{x}_{ti}\mathbf{z}_i
      -\frac{1}{\gamma_0}\mathbf{x}'_t\boldsymbol{\delta}_0\mathbf{z}_{K+1}\biggr].$

Define $g(\mathbf{x}_t,\mathbf{z},\boldsymbol{\delta}_0,\gamma_0)$
as
$\displaystyle\biggl[-\sum_{i=1}^{K}\lambda(v_{0t})\mathbf{x}_{ti}\mathbf{z}_i
      +\bigl(\frac{1}{\gamma_0}+\lambda(v_{0t})c\bigr)\mathbf{z}_{K+1}\biggr]$.

Then (2.25) is equal to
\begin{equation}
\gamma_0\epsilon_tf(\mathbf{x}_t,\mathbf{z},\boldsymbol{\delta}_0,\gamma_0)
 1_{\{\epsilon_t>c-\frac{1}{\gamma_0}\mathbf{x}'_t\boldsymbol{\delta}_0\}}
      -\gamma_0\epsilon_t^2\mathbf{z}_{K+1}1_{\{\epsilon_t>c-\frac{1}{\gamma_0}\mathbf{x}'_t\boldsymbol{\delta}_0\}}
      +g(\mathbf{x}_t,\mathbf{z},\boldsymbol{\delta}_0,\gamma_0)1_{\{\epsilon_t>c-\frac{1}{\gamma_0}\mathbf{x}'_t\boldsymbol{\delta}_0\}}.
\end{equation}
Thus (2.24) is equivalent to
\begin{equation}
\gamma_0\epsilon_tf(\mathbf{x}_t,\mathbf{z},\boldsymbol{\delta}_0,\gamma_0)
 1_{\{\epsilon_t>c-\frac{1}{\gamma_0}\mathbf{x}'_t\boldsymbol{\delta}_0\}}
      -\gamma_0\epsilon_t^2\mathbf{z}_{K+1}1_{\{\epsilon_t>c-\frac{1}{\gamma_0}\mathbf{x}'_t\boldsymbol{\delta}_0\}}
      +g(\mathbf{x}_t,\mathbf{z},\boldsymbol{\delta}_0,\gamma_0)1_{\{\epsilon_t>c-\frac{1}{\gamma_0}\mathbf{x}'_t\boldsymbol{\delta}_0\}}=0\,\,\,\textrm{a.e.},
\end{equation}
namely,
\begin{equation}
-\epsilon_t^2\mathbf{z}_{K+1}1_{\{\epsilon_t>c-\frac{1}{\gamma_0}\mathbf{x}'_t\boldsymbol{\delta}_0\}}
+\epsilon_tf(\mathbf{x}_t,\mathbf{z},\boldsymbol{\delta}_0,\gamma_0)
 1_{\{\epsilon_t>c-\frac{1}{\gamma_0}\mathbf{x}'_t\boldsymbol{\delta}_0\}}
+\frac{1}{\gamma_0}g(\mathbf{x}_t,\mathbf{z},\boldsymbol{\delta}_0,\gamma_0)1_{\{\epsilon_t>c-\frac{1}{\gamma_0}\mathbf{x}'_t\boldsymbol{\delta}_0\}}=0\,\,\,\textrm{a.e.}.
\end{equation}

Suppose $\mathbf{z}_{K+1}\neq 0$. When
$\epsilon_t>c-\frac{1}{\gamma_0}\mathbf{x}'_t\boldsymbol{\delta}_0$
holds, (2.28) is a quadratic equation of $\epsilon_t$. If the
quadratic equation has solutions, then the solutions are
\begin{equation}
\epsilon_t=\frac{-f(\mathbf{x}_t,\mathbf{z},\boldsymbol{\delta}_0,\gamma_0)\pm\sqrt{f^2(\mathbf{x}_t,\mathbf{z},\boldsymbol{\delta}_0,\gamma_0)
+\frac{4\mathbf{z}_{K+1}}{\gamma_0}g(\mathbf{x}_t,\mathbf{z},\boldsymbol{\delta}_0,\gamma_0)}}{-2\mathbf{z}_{K+1}}\,\,\,\,\,\,\textrm{a.e.}.
\end{equation}
But this contradicts the fact that $\epsilon_t|\mathbf{x}_t\sim
N(0,\sigma^2_{0})$ (see (1.2)).

Therefore $\mathbf{z}_{K+1}$ must be $0$. This means that
\begin{equation}
{\mathbf{z}}'E[\mathbf{H}(\mathbf{w}_t;\boldsymbol{\delta}_0,
\gamma_0)]\mathbf{z}=0\,\, \Longrightarrow\,\,
\mathbf{z}=(\mathbf{z}_1,...,\mathbf{z}_K,0)'.
\end{equation}

By (1.9),
\begin{align}
&(\mathbf{z}_1,...,\mathbf{z}_K,0)E[\mathbf{H}(\mathbf{w}_t;\boldsymbol{\delta}_0,
\gamma_0)](\mathbf{z}_1,...,\mathbf{z}_K,0)'\notag\\
=&-(\mathbf{z}_1,...,\mathbf{z}_K)
E\Bigl[\bigl(1-\lambda(v_{0t})[\lambda(v_{0t})-v_{0t}]\bigr)\mathbf{x_tx'_t}\Bigr](\mathbf{z}_1,...,\mathbf{z}_K)'\notag\\
=&-E\Bigl[(\mathbf{z}_1,...,\mathbf{z}_K)
\bigl(1-\lambda(v_{0t})[\lambda(v_{0t})-v_{0t}]\bigr)\mathbf{x_tx'_t}(\mathbf{z}_1,...,\mathbf{z}_K)'\Bigr],
\end{align}
We know that for any positive constant $\overline{v} \in \R$,
\begin{align}
&\,E\Bigl[(\mathbf{z}_1,...,\mathbf{z}_K)
\bigl(1-\lambda(v_{0t})[\lambda(v_{0t})-v_{0t}]\bigr)\mathbf{x_tx'_t}(\mathbf{z}_1,...,\mathbf{z}_K)'\Bigr]\notag\\
\geq&\,E\Bigl[(\mathbf{z}_1,...,\mathbf{z}_K) 1_{\{|v_{0t}|\leq
{\overline{v}}\}}\bigl(1-\lambda(v_{0t})[\lambda(v_{0t})-v_{0t}]\bigr)\mathbf{x_tx'_t}(\mathbf{z}_1,...,\mathbf{z}_K)'\Bigr]\notag\\
\geq&\,\overline{C}E\Bigl[(\mathbf{z}_1,...,\mathbf{z}_K)
1_{\{|v_{0t}|\leq
{\overline{v}}\}}\mathbf{x_tx'_t}(\mathbf{z}_1,...,\mathbf{z}_K)'\Bigr],
\end{align}
where $\overline{C}$ is a positive constant depending on
$\overline{v}$. Here we used the fact that
$\lambda(v_{0t})[\lambda(v_{0t})-v_{0t}]$ is between $0$ and $1$,
and asymptotes to $0$ as $v_{0t}\rightarrow -\infty$ and to $1$ as
$v_{0t}\rightarrow \infty$.

Therefore, by (2.30), (2.31) and (2.32),
\begin{equation}
{\mathbf{z}}'E[\mathbf{H}(\mathbf{w}_t;\boldsymbol{\delta}_0,
\gamma_0)]\mathbf{z}=0\,\, \Longrightarrow\,\,
(\mathbf{z}_1,...,\mathbf{z}_K)E[ 1_{\{|v_{0t}|\leq
{\overline{v}}\}}\mathbf{x_tx'_t}](\mathbf{z}_1,...,\mathbf{z}_K)'=0.
\end{equation}

It is easy to see that if $E[\mathbf{x}_t\mathbf{x}'_t]$ is
nonsingular, then for large enough $\overline{v}$,
$E[1_{\{|v_{0t}|\leq {\overline{v}}\}}\mathbf{x}_t\mathbf{x}'_t]$
is also nonsingular. Thus when $E[\mathbf{x}_t\mathbf{x}'_t]$ is
nonsingular, for large enough $\overline{v}$,
$\{\mathbf{z}_i=0,\,i=1,...,K\}$ is the only solution satisfying
\begin{equation}
(\mathbf{z}_1,...,\mathbf{z}_K)E[ 1_{\{|v_{0t}|\leq
{\overline{v}}\}}\mathbf{x}_t\mathbf{x}'_t](\mathbf{z}_1,...,\mathbf{z}_K)'=0.
\end{equation}

This means that when $E[\mathbf{x}_t\mathbf{x}'_t]$ is
nonsingular,
\begin{equation}
{\mathbf{z}}'E[\mathbf{H}(\mathbf{w}_t;\boldsymbol{\delta}_0,
\gamma_0)]\mathbf{z}=0\,\, \Longrightarrow\,\,
\mathbf{z}=(0,...,0,0)'.
\end{equation}
Thus when $E[\mathbf{x_tx'_t}]$ is nonsingular,
$E[\mathbf{H}(\mathbf{w}_t;\boldsymbol{\delta}_0, \gamma_0)]$ is
nonsingular, i.e. condition 5 of Proposition 1.3 holds. \qed

Theorem 2.1 implies that when $E[\mathbf{x_tx'_t}]$ is
nonsingular, under Assumptions $1$ and $2$, the conditional ML
estimator $(\boldsymbol{\hat{\delta}},\hat{\gamma})$ is asymptotic
normal with
$\texttt{Avar}(\boldsymbol{\hat{\delta}},\hat{\gamma})$ given by
the following:
\begin{equation}
\texttt{Avar}(\boldsymbol{\hat{\delta}},\hat{\gamma})=-\{E[\mathbf{H}(\mathbf{w}_t;\boldsymbol{\delta}_0,\gamma_0)]\}^{-1}=\{E[\mathbf{s}(\mathbf{w}_t;\boldsymbol{\delta}_0,\gamma_0)
\mathbf{s}(\mathbf{w}_t;\boldsymbol{\delta}_0,\gamma_0)']\}^{-1}.
\end{equation}

To recover original parameters and obtain the asymptotic variance
of $(\boldsymbol{\hat{\beta}},\hat{\sigma}^2)$, the delta method
can be applied. (see page 517 of Hayashi (2000))

\section{Asymptotic Normality of the Conditional ML of the
Tobit Model}

In this section, we show that for the Tobit Model introduced in
subsection 1.2, the conditions of Proposition 1.3 are satisfied,
thereafter the asymptotic normality of the conditional ML is
verified.

It is easy to see that conditions $1$ and $2$ of Proposition 1.3
are satisfied.

We know that in the expression of
$\mathbf{H}(\mathbf{w}_t;\boldsymbol{\delta },\gamma)$ of (1.18),
$1-D_t$ is either $1$ or $0$, so bounded, and
$\lambda(-v_t)[\lambda(-v_t)+v_t]$ is between $0$ and $1$, so
bounded. Thus by the same argument as in (2.1) through (2.18) for
the truncated regression model, we can show that condition $4$ of
Proposition 1.3 is satisfied.

In the following, we show conditions $3$ and $5$ of Proposition
1.3 are satisfied.

Before we move on to the main results,  we need some preliminary
results on the conditional moments of
$\gamma_0y_t-\mathbf{x}'_t\boldsymbol{\delta}_0$ which is
conditioning on $\{\mathbf{x}_t,y_t>c\}$.
\begin{lemma}
For the Tobit model, the following equalities hold:
\begin{align}
E[\gamma_0y_t-\mathbf{x}'_t\boldsymbol{\delta}_0|\mathbf{x}_t,y_t>c]&=\lambda(v_{0t}),\\
E\bigl[(\gamma_0y_t-\mathbf{x}'_t\boldsymbol{\delta}_0)^2\mid\mathbf{x}_t,y_t>c\bigr]&=v_{0t}\lambda(v_{0t})+1,\\
E\bigl[(\gamma_0y_t-\mathbf{x}'_t\boldsymbol{\delta}_0)^3|\mathbf{x}_t,y_t>c\bigr]&=v^2_{0t}\lambda(v_{0t})+2\lambda(v_{0t}),\\
E\bigl[(\gamma_0y_t-\mathbf{x}'_t\boldsymbol{\delta}_0)^4|\mathbf{x}_t,y_t>c\bigr]&=v^3_{0t}\lambda(v_{0t})+3[v_{0t}\lambda(v_{0t})+1].
\end{align}
\end{lemma}
\pf

\begin{align}
E[\gamma_0y_t-\mathbf{x}'_t\boldsymbol{\delta}_0|\mathbf{x}_t,y_t>c]
&=E\Bigl[\frac{y_t-\mathbf{x}'_t\boldsymbol{\beta}_0}{\sigma_0}\mid\mathbf{x}_t,\,\,
\tfrac{y_t-\mathbf{x}'_t\boldsymbol{\beta}_0}{\sigma_0}>\tfrac{c-\mathbf{x}'_t\boldsymbol{\beta}_0}{\sigma_0}       \Bigr]\notag\\
&=\int_{\frac{c-\mathbf{x}'_t\boldsymbol{\beta}_0}{\sigma_0}}
   ^{\infty}\tilde{y}\,\frac{\phi(\tilde{y})}{1-\Phi(\frac{c-\mathbf{x}'_t\boldsymbol{\beta}_0}{\sigma_0})}\,d\tilde{y}\notag\\
&=\int_{\frac{c-\mathbf{x}'_t\boldsymbol{\beta}_0}{\sigma_0}}
   ^{\infty}\tilde{y}\,\frac{\frac{1}{\sqrt{2\pi}}e^{-\frac{\tilde{y}^2}{2}}}{1-\Phi(\frac{c-\mathbf{x}'_t\boldsymbol{\beta}_0}{\sigma_0})}\,d\tilde{y}\notag\\
&=-\int_{\frac{c-\mathbf{x}'_t\boldsymbol{\beta}_0}{\sigma_0}}
   ^{\infty}\,\frac{\frac{1}{\sqrt{2\pi}}}{1-\Phi(\frac{c-\mathbf{x}'_t\boldsymbol{\beta}_0}{\sigma_0})}\,de^{-\frac{\tilde{y}^2}{2}}\notag\\
&=\frac{\phi(\frac{c-\mathbf{x}'_t\boldsymbol{\beta}_0}{\sigma_0})}{1-\Phi(\frac{c-\mathbf{x}'_t\boldsymbol{\beta}_0}{\sigma_0})}\notag\\
&=\lambda(v_{0t}).\notag
\end{align}
This shows that (3.1) holds.

\begin{align}
E\bigl[(\gamma_0y_t-\mathbf{x}'_t\boldsymbol{\delta}_0)^2\mid\mathbf{x}_t,y_t>c\bigr]
&=E\biggl[\Bigl(\dfrac{y_t-\mathbf{x}'_t\boldsymbol{\beta}_0}{\sigma_0}\Bigr)^2\mid\mathbf{x}_t,\,\,
\tfrac{y_t-\mathbf{x}'_t\boldsymbol{\beta}_0}{\sigma_0}>\tfrac{c-\mathbf{x}'_t\boldsymbol{\beta}_0}{\sigma_0}\biggr]\notag\\
&=\int_{\frac{c-\mathbf{x}'_t\boldsymbol{\beta}_0}{\sigma_0}}
   ^{\infty}\tilde{y}^2\,\frac{\phi(\tilde{y})}{1-\Phi(\frac{c-\mathbf{x}'_t\boldsymbol{\beta}_0}{\sigma_0})}\,d\tilde{y}\notag\\
&=\int_{\frac{c-\mathbf{x}'_t\boldsymbol{\beta}_0}{\sigma_0}}
   ^{\infty}\tilde{y}^2\,\frac{\frac{1}{\sqrt{2\pi}}e^{-\frac{\tilde{y}^2}{2}}}{1-\Phi(\frac{c-\mathbf{x}'_t\boldsymbol{\beta}_0}{\sigma_0})}\,d\tilde{y}\notag\\
      &=-\int_{\frac{c-\mathbf{x}'_t\boldsymbol{\beta}_0}{\sigma_0}}
   ^{\infty}\,\frac{\tilde{y}\,\frac{1}{\sqrt{2\pi}}}{1-\Phi(\frac{c-\mathbf{x}'_t\boldsymbol{\beta}_0}{\sigma_0})}\,de^{-\frac{\tilde{y}^2}{2}}\notag\\
&=\Bigl(\dfrac{c-\mathbf{x}'_t\boldsymbol{\beta}_0}{\sigma_0}\Bigr)
\frac{\phi(\frac{c-\mathbf{x}'_t\boldsymbol{\beta}_0}{\sigma_0})}{1-\Phi(\frac{c-\mathbf{x}'_t\boldsymbol{\beta}_0}{\sigma_0})}
+\int_{\frac{c-\mathbf{x}'_t\boldsymbol{\beta}_0}{\sigma_0}}
   ^{\infty}\frac{\frac{1}{\sqrt{2\pi}}e^{-\frac{\tilde{y}^2}{2}}}{1-\Phi(\frac{c-\mathbf{x}'_t\boldsymbol{\beta}_0}{\sigma_0})}\,d\tilde{y}\notag\\
&=v_{0t}\lambda(v_{0t})+1.\notag
\end{align}
Thus (3.2) holds.

\begin{align}
E\bigl[(\gamma_0y_t-\mathbf{x}'_t\boldsymbol{\delta}_0)^3|\mathbf{x}_t,y_t>c\bigr]
&=E\biggl[\Bigl(\dfrac{y_t-\mathbf{x}'_t\boldsymbol{\beta}_0}{\sigma_0}\Bigr)^3\mid\mathbf{x}_t,\,\,
\tfrac{y_t-\mathbf{x}'_t\boldsymbol{\beta}_0}{\sigma_0}>\tfrac{c-\mathbf{x}'_t\boldsymbol{\beta}_0}{\sigma_0}\biggr]\notag\\
&=\int_{\frac{c-\mathbf{x}'_t\boldsymbol{\beta}_0}{\sigma_0}}
   ^{\infty}\tilde{y}^3\,\frac{\phi(\tilde{y})}{1-\Phi(\frac{c-\mathbf{x}'_t\boldsymbol{\beta}_0}{\sigma_0})}\,d\tilde{y}\notag\\
&=\int_{\frac{c-\mathbf{x}'_t\boldsymbol{\beta}_0}{\sigma_0}}
   ^{\infty}\tilde{y}^3\,\frac{\frac{1}{\sqrt{2\pi}}e^{-\frac{\tilde{y}^2}{2}}}{1-\Phi(\frac{c-\mathbf{x}'_t\boldsymbol{\beta}_0}{\sigma_0})}\,d\tilde{y}\notag\\
&=-\int_{\frac{c-\mathbf{x}'_t\boldsymbol{\beta}_0}{\sigma_0}}
   ^{\infty}\,\frac{\tilde{y}^2\,\frac{1}{\sqrt{2\pi}}}{1-\Phi(\frac{c-\mathbf{x}'_t\boldsymbol{\beta}_0}{\sigma_0})}\,de^{-\frac{\tilde{y}^2}{2}}\notag\\
&=\Bigl(\dfrac{c-\mathbf{x}'_t\boldsymbol{\beta}_0}{\sigma_0}\Bigr)^2
\frac{\phi(\frac{c-\mathbf{x}'_t\boldsymbol{\beta}_0}{\sigma_0})}{1-\Phi(\frac{c-\mathbf{x}'_t\boldsymbol{\beta}_0}{\sigma_0})}
+2\int_{\frac{c-\mathbf{x}'_t\boldsymbol{\beta}_0}{\sigma_0}}
   ^{\infty}\tilde{y}\frac{\frac{1}{\sqrt{2\pi}}e^{-\frac{\tilde{y}^2}{2}}}{1-\Phi(\frac{c-\mathbf{x}'_t\boldsymbol{\beta}_0}{\sigma_0})}\,d\tilde{y}\notag\\
&=v^2_{0t}\lambda(v_{0t})+2\lambda(v_{0t}),\textrm{ }\textrm{
}\qquad\textrm{ ( by (3.1) )}\notag
\end{align}
This shows that (3.3) holds.

\begin{align}
&E\bigl[(\gamma_0y_t-\mathbf{x}'_t\boldsymbol{\delta}_0)^4|\mathbf{x}_t,y_t>c
\bigr]\notag\\
&=E\biggl[\Bigl(\dfrac{y_t-\mathbf{x}'_t\boldsymbol{\beta}_0}{\sigma_0}\Bigr)^4\mid\mathbf{x}_t,\,\,
\tfrac{y_t-\mathbf{x}'_t\boldsymbol{\beta}_0}{\sigma_0}>\tfrac{c-\mathbf{x}'_t\boldsymbol{\beta}_0}{\sigma_0}\biggr]\notag\\
&=\int_{\frac{c-\mathbf{x}'_t\boldsymbol{\beta}_0}{\sigma_0}}
   ^{\infty}\tilde{y}^4\,\frac{\phi(\tilde{y})}{1-\Phi(\frac{c-\mathbf{x}'_t\boldsymbol{\beta}_0}{\sigma_0})}\,d\tilde{y}\notag\\
&=\int_{\frac{c-\mathbf{x}'_t\boldsymbol{\beta}_0}{\sigma_0}}
   ^{\infty}\tilde{y}^4\,\frac{\frac{1}{\sqrt{2\pi}}e^{-\frac{\tilde{y}^2}{2}}}{1-\Phi(\frac{c-\mathbf{x}'_t\boldsymbol{\beta}_0}{\sigma_0})}\,d\tilde{y}\notag\\
&=-\int_{\frac{c-\mathbf{x}'_t\boldsymbol{\beta}_0}{\sigma_0}}
   ^{\infty}\,\frac{\tilde{y}^3\,\frac{1}{\sqrt{2\pi}}}{1-\Phi(\frac{c-\mathbf{x}'_t\boldsymbol{\beta}_0}{\sigma_0})}\,de^{-\frac{\tilde{y}^2}{2}}\notag\\
&=\Bigl(\dfrac{c-\mathbf{x}'_t\boldsymbol{\beta}_0}{\sigma_0}\Bigr)^3
\frac{\phi(\frac{c-\mathbf{x}'_t\boldsymbol{\beta}_0}{\sigma_0})}{1-\Phi(\frac{c-\mathbf{x}'_t\boldsymbol{\beta}_0}{\sigma_0})}
+3\int_{\frac{c-\mathbf{x}'_t\boldsymbol{\beta}_0}{\sigma_0}}
   ^{\infty}\tilde{y}^2\frac{\frac{1}{\sqrt{2\pi}}e^{-\frac{\tilde{y}^2}{2}}}{1-\Phi(\frac{c-\mathbf{x}'_t\boldsymbol{\beta}_0}{\sigma_0})}\,d\tilde{y}\notag\\
&=v^3_{0t}\lambda(v_{0t})+3[v_{0t}\lambda(v_{0t})+1],\textrm{
}\textrm{ }\qquad\textrm{ ( by (3.2) )}\notag
\end{align}
Thus (3.4) holds. \qed

Next we claim
\begin{thm}
For the Tobit Model satisfying Assumption 1', if $\mathbf{
E(}\mathbf{x}_t\mathbf{x}'_t\mathbf{)}$ is nonsingular, then
condition $3$ of Proposition 1.3 is satisfied.
\end{thm}
\pf

First we show that
$E[\mathbf{s}(\mathbf{w}_t;\boldsymbol{\delta}_0,\gamma_0)|
\mathbf{x}_t]=0$.

In terms of the expression of
$\mathbf{s}(\mathbf{w}_t;\boldsymbol{\delta},\gamma)$ in (1.17),
\begin{equation}
E[\mathbf{s}(\mathbf{w}_t;\boldsymbol{\delta}_0,\gamma_0)|
\mathbf{x}_t]=\begin{bmatrix}E\bigl[(1-D_t)(\gamma_0
y_t-\mathbf{x}'_t\boldsymbol{\delta}_0)\mathbf{x}_t+D_t\lambda(-v_{0t})(-\mathbf{x}_t)|\mathbf{x}_t\bigr]\\
\textrm{ }\\
 E\Bigl[(1-D_t)\bigl[\tfrac{1}{\gamma_0}-(\gamma_0
y_t-\mathbf{x}'_t\boldsymbol{\delta}_0)y_t\bigr]+D_t\lambda(-v_{0t})c|\mathbf{x}_t\Bigr]
\end{bmatrix},
\end{equation}
where
\begin{align}
&E\bigl[(1-D_t)(\gamma_0
y_t-\mathbf{x}'_t\boldsymbol{\delta}_0)\mathbf{x}_t+D_t\lambda(-v_{0t})(-\mathbf{x}_t)|\mathbf{x}_t\bigr]\notag\\
&=E\bigl[(\gamma_0
y_t-\mathbf{x}'_t\boldsymbol{\delta}_0)\mathbf{x}_t|\mathbf{x}_t,y_t>c\bigr]\,\texttt{Prob}[y_t>c|\mathbf{x}_t]
+E\bigl[\lambda(-v_{0t})(-\mathbf{x}_t)|\mathbf{x}_t,y_t=c\bigr]\,\texttt{Prob}[y_t= c|\mathbf{x}_t]\notag\\
&=E\bigl[(\gamma_0
y_t-\mathbf{x}'_t\boldsymbol{\delta}_0)|\mathbf{x}_t,y_t>c\bigr](\mathbf{x}_t)\,\texttt{Prob}[y_t>c|\mathbf{x}_t]
-\lambda(-v_{0t})(\mathbf{x}_t)\,\texttt{Prob}[y_t=
c|\mathbf{x}_t]\notag
\end{align}
\begin{align}
&=\lambda(v_{0t})(\mathbf{x}_t)\,\texttt{Prob}[y_t>c|\mathbf{x}_t]-\lambda(-v_{0t})(\mathbf{x}_t)\,\texttt{Prob}[y_t= c|\mathbf{x}_t]\qquad \qquad \textrm{ ( by (3.1) )}\notag\\
&=\tfrac{\phi(\frac{c-\mathbf{x}'_t\boldsymbol{\beta}_0}{\sigma_0})}{1-\Phi(\tfrac{c-\mathbf{x}'_t\boldsymbol{\beta}_0}{\sigma_0})}
(\mathbf{x}_t)\Bigl[1-\Phi\Bigl(\tfrac{c-\mathbf{x}'_t\boldsymbol{\beta}_0}{\sigma_0}\Bigr)\Bigr]-
\tfrac{\phi(-\frac{c-\mathbf{x}'_t\boldsymbol{\beta}_0}{\sigma_0})}{1-\Phi(-\tfrac{c-\mathbf{x}'_t\boldsymbol{\beta}_0}{\sigma_0})}
(\mathbf{x}_t)\Phi\Bigl(\tfrac{c-\mathbf{x}'_t\boldsymbol{\beta}_0}{\sigma_0}\Bigr)\notag\\
&=\phi\Big(\tfrac{c-\mathbf{x}'_t\boldsymbol{\beta}_0}{\sigma_0}\Bigr)(\mathbf{x}_t)-\phi\Bigl(-\tfrac{c-\mathbf{x}'_t\boldsymbol{\beta}_0}{\sigma_0}\Bigr)
(\mathbf{x}_t)\notag\\
&=0,
\end{align}

and
\begin{align}
&E\Bigl[(1-D_t)\bigl[\tfrac{1}{\gamma_0}-(\gamma_0
y_t-\mathbf{x}'_t\boldsymbol{\delta}_0)y_t\bigr]+D_t\lambda(-v_{0t})c\mid\mathbf{x}_t\Bigr]\notag\\
&=E\bigl[\tfrac{1}{\gamma_0}-(\gamma_0
y_t-\mathbf{x}'_t\boldsymbol{\delta}_0)y_t|\mathbf{x}_t,y_t>c\bigr]\,\texttt{Prob}[y_t>c|\mathbf{x}_t]
+E\bigl[\lambda(-v_{0t})c|\mathbf{x}_t,y_t=c\bigr]\,\texttt{Prob}[y_t= c|\mathbf{x}_t]\notag\\
&=E\bigl[\tfrac{1}{\gamma_0}-(\gamma_0
y_t-\mathbf{x}'_t\boldsymbol{\delta}_0)\tfrac{1}{\gamma_0}(\gamma_0y_t-\mathbf{x}'_t\boldsymbol{\delta}_0
+\mathbf{x}'_t\boldsymbol{\delta}_0)|\mathbf{x}_t,y_t>c\bigr]\,\texttt{Prob}[y_t>c|\mathbf{x}_t]\notag\\
&\quad+E\bigl[\lambda(-v_{0t})c|\mathbf{x}_t,y_t=c\bigr]\,\texttt{Prob}[y_t= c|\mathbf{x}_t]\notag\\
&=\Bigl\{\tfrac{1}{\gamma_0}-\tfrac{1}{\gamma_0}E\bigl[(\gamma_0
y_t-\mathbf{x}'_t\boldsymbol{\delta}_0)^2|\mathbf{x}_t,y_t>c\bigr]-\tfrac{1}{\gamma_0}
E\bigl[(\gamma_0
y_t-\mathbf{x}'_t\boldsymbol{\delta}_0)|\mathbf{x}_t,y_t>c\bigr](\mathbf{x}'_t\boldsymbol{\delta}_0)\Bigr\}\,\notag\\
&\quad \cdot \texttt{Prob}[y_t>c|\mathbf{x}_t]+\lambda(-v_{0t})c\,\texttt{Prob}[y_t= c|\mathbf{x}_t]\notag\\
&=\Bigl\{\tfrac{1}{\gamma_0}-\tfrac{1}{\gamma_0}[v_{0t}\lambda(v_{0t})+1]-\tfrac{1}{\gamma_0}
\lambda(v_{0t})(\mathbf{x}'_t\boldsymbol{\delta}_0)\Bigr\}\,\texttt{Prob}[y_t>c|\mathbf{x}_t]\qquad \qquad \textrm{ ( by (3.2) and (3.1) )}\notag\\
&\quad+\lambda(-v_{0t})c\,\texttt{Prob}[y_t=c|\mathbf{x}_t]\notag\\
&=-\lambda(v_{0t})c\,\texttt{Prob}[y_t>c|\mathbf{x}_t]+\lambda(-v_{0t})c\,\texttt{Prob}[y_t=
c|\mathbf{x}_t]\quad\quad\quad \,\, (\textrm{ by
}v_{0t}+\mathbf{x}'_t\boldsymbol{\delta}_0=\gamma_0c \textrm{ })\notag\\
&=-\tfrac{\phi(\frac{c-\mathbf{x}'_t\boldsymbol{\beta}_0}{\sigma_0})}{1-\Phi(\tfrac{c-\mathbf{x}'_t\boldsymbol{\beta}_0}{\sigma_0})}
c\Bigl[1-\Phi\Bigl(\tfrac{c-\mathbf{x}'_t\boldsymbol{\beta}_0}{\sigma_0}\Bigr)\Bigr]+
\tfrac{\phi(-\frac{c-\mathbf{x}'_t\boldsymbol{\beta}_0}{\sigma_0})}{1-\Phi(-\tfrac{c-\mathbf{x}'_t\boldsymbol{\beta}_0}{\sigma_0})}
c\,\Phi\Bigl(\tfrac{c-\mathbf{x}'_t\boldsymbol{\beta}_0}{\sigma_0}\Bigr)\notag\\
&=-\phi\Big(\tfrac{c-\mathbf{x}'_t\boldsymbol{\beta}_0}{\sigma_0}\Bigr)c+\phi\Bigl(-\tfrac{c-\mathbf{x}'_t\boldsymbol{\beta}_0}{\sigma_0}\Bigr)
c\notag\\
&=0.
\end{align}
Therefore
$E[\mathbf{s}(\mathbf{w}_t;\boldsymbol{\delta}_0,\gamma_0)|
\mathbf{x}_t]=0$.

Next we claim
$E[\mathbf{s}(\mathbf{w}_t;\boldsymbol{\delta}_0,\gamma_0)\mathbf{s}(\mathbf{w}_t;\boldsymbol{\delta}_0,\gamma_0)'|
\mathbf{x}_t]=-E[\mathbf{H}(\mathbf{w}_t;\boldsymbol{\delta}_0,\gamma_0)|\mathbf{x}_t]$.

In terms of the expression of
$\mathbf{s}(\mathbf{w}_t;\boldsymbol{\delta},\gamma)$ in (1.17),
we write the matrix\\
$\mathbf{s}(\mathbf{w}_t;\boldsymbol{\delta}_0,\gamma_0)\mathbf{s}(\mathbf{w}_t;\boldsymbol{\delta}_0,\gamma_0)'$
as
\begin{equation}
\begin{bmatrix}
\mathbf{M}_{11}&\mathbf{M}_{12}\\
\mathbf{M}_{21}&\mathbf{M}_{22}
\end{bmatrix},
\end{equation}
where
\begin{align}
\mathbf{M}_{11}&=(1-D_t)(\gamma_0y_t-\mathbf{x}'_t\mathbf{\delta}_0)^2\mathbf{x_tx'_t}+D_t\lambda^2(-v_{0t})\mathbf{x_tx'_t},\\
\mathbf{M}_{12}&=(1-D_t)(\gamma_0y_t-\mathbf{x}'_t\mathbf{\delta}_0)\mathbf{x}_t
\bigl[\tfrac{1}{\gamma_0}-(\gamma_0y_t-\mathbf{x}'_t\mathbf{\delta}_0)y_t\bigr]
+D_t\lambda^2(-v_{0t})(-\mathbf{x}_tc),\\
\mathbf{M}_{21}&=\mathbf{M}'_{12},\\
\mathbf{M}_{22}&=(1-D_t)\bigl[\tfrac{1}{\gamma_0}-(\gamma_0y_t-\mathbf{x}'_t\mathbf{\delta}_0)y_t\bigr]^2
+D_t\lambda^2(-v_{0t})c^2.
\end{align}

In terms of the expression of
$\mathbf{H}(\mathbf{w}_t;\boldsymbol{\delta},\gamma)$ in (1.18),
we write the matrix
$-\mathbf{H}(\mathbf{w}_t;\boldsymbol{\delta}_0,\gamma_0)$ as
\begin{equation}
\begin{bmatrix}
\mathbf{N}_{11}&\mathbf{N}_{12}\\
\mathbf{N}_{21}&\mathbf{N}_{22}
\end{bmatrix},
\end{equation}
where
\begin{align}
\mathbf{N}_{11}&=(1-D_t)\mathbf{x_tx'_t}+D_t\lambda(-v_{0t})[\lambda(-v_{0t})+v_{0t}]\mathbf{x_tx'_t},\\
\mathbf{N}_{12}&=(1-D_t)(-y_t\mathbf{x}_t)+D_t\lambda(-v_{0t})[\lambda(-v_{0t})+v_{0t}](-c\mathbf{x}_t),\\
\mathbf{N}_{21}&=\mathbf{N}'_{21},\\
\mathbf{N}_{22}&=(1-D_t)\bigl[\tfrac{1}{\gamma^2_0}+y^2_t\bigr]+D_t\lambda(-v_{0t})[\lambda(-v_{0t})+v_{0t}]c^2.
\end{align}

If we can show
$E[\mathbf{N}_{ij}|\mathbf{x}_t]=E[\mathbf{M}_{ij}|\mathbf{x}_t],\,\,\,i,=1,2,\,j=1,2$,\,
then the claim is true. The verification is as follows.

Since
\begin{align}
E[\mathbf{M}_{11}|\mathbf{x}_t]&=E\bigl[(\gamma_0y_t-\mathbf{x}'_t\mathbf{\delta}_0)^2|\mathbf{x}_t,y_t>c\bigr]\mathbf{x_tx'_t}\,\texttt{Prob}[y_t>c|\mathbf{x}_t]
+\lambda^2(-v_{0t})\mathbf{x_tx'_t}\,\texttt{Prob}[y_t=c|\mathbf{x}_t]\notag\\
&=\Bigl\{\{[v_{0t}\lambda(v_{0t})+1]\Bigl[1-\Phi\bigl(\tfrac{c-\mathbf{x}'_t\boldsymbol{\beta}_0}{\sigma_0}\bigr)\Bigr]
+\lambda^2(-v_{0t})\Phi\bigl(\tfrac{c-\mathbf{x}'_t\boldsymbol{\beta}_0}{\sigma_0}\bigr)\Bigr\}\mathbf{x_tx'_t}
\qquad \quad \textrm { ( by (3.2) )}\notag\\
&=\Bigl\{v_{0t}\lambda(v_{0t})\Bigl[1-\Phi\bigl(\tfrac{c-\mathbf{x}'_t\boldsymbol{\beta}_0}{\sigma_0}\bigr)\Bigr]
+\Bigl[1-\Phi\bigl(\tfrac{c-\mathbf{x}'_t\boldsymbol{\beta}_0}{\sigma_0}\bigr)\Bigr]
+\lambda^2(-v_{0t})\Phi\bigl(\tfrac{c-\mathbf{x}'_t\boldsymbol{\beta}_0}{\sigma_0}\bigr)\Bigr\}\mathbf{x_tx'_t}\notag\\
&=\Bigl\{v_{0t}\phi
\bigl(\tfrac{c-\mathbf{x}'_t\boldsymbol{\beta}_0}{\sigma_0}\bigr)+\Bigl[1-\Phi\bigl(\tfrac{c-\mathbf{x}'_t\boldsymbol{\beta}_0}{\sigma_0}\bigr)\Bigr]
+\lambda^2(-v_{0t})\Phi\bigl(\tfrac{c-\mathbf{x}'_t\boldsymbol{\beta}_0}{\sigma_0}\bigr)\Bigr\}\mathbf{x_tx'_t},
\end{align}
and
\begin{align}
E[\mathbf{N}_{11}|\mathbf{x}_t]&=\mathbf{x_tx'_t}\,\texttt{Prob}[y_t>c|\mathbf{x}_t]
+\lambda(-v_{0t})[\lambda(-v_{0t})+v_{0t}]\mathbf{x_tx'_t}\,\texttt{Prob}[y_t=c|\mathbf{x}_t]\notag\\
&=\Bigl\{\Bigl[1-\Phi\bigl(\tfrac{c-\mathbf{x}'_t\boldsymbol{\beta}_0}{\sigma_0}\bigr)\Bigr]
+[\lambda^2(-v_{0t})+\lambda(-v_{0t})v_{0t}]\Phi\bigl(\tfrac{c-\mathbf{x}'_t\boldsymbol{\beta}_0}{\sigma_0}\bigr)\Bigr\}\mathbf{x_tx'_t}\notag\\
&=\Bigl\{\Bigl[1-\Phi\bigl(\tfrac{c-\mathbf{x}'_t\boldsymbol{\beta}_0}{\sigma_0}\bigr)\Bigr]
+\lambda^2(-v_{0t})\Phi\bigl(\tfrac{c-\mathbf{x}'_t\boldsymbol{\beta}_0}{\sigma_0}\bigr)+\lambda(-v_{0t})v_{0t}
\Phi\bigl(\tfrac{c-\mathbf{x}'_t\boldsymbol{\beta}_0}{\sigma_0}\bigr)\Bigr\}\mathbf{x_tx'_t}\notag\\
&=\Bigl\{\Bigl[1-\Phi\bigl(\tfrac{c-\mathbf{x}'_t\boldsymbol{\beta}_0}{\sigma_0}\bigr)\Bigr]
+\lambda^2(-v_{0t})\Phi\bigl(\tfrac{c-\mathbf{x}'_t\boldsymbol{\beta}_0}{\sigma_0}\bigr)+v_{0t}
\phi\bigl(\tfrac{c-\mathbf{x}'_t\boldsymbol{\beta}_0}{\sigma_0}\bigr)\Bigr\}\mathbf{x_tx'_t},\notag\\
&\qquad (\textrm{ by the fact that
}\lambda(-v_{0t})\Phi\bigl(\tfrac{c-\mathbf{x}'_t\boldsymbol{\beta}_0}{\sigma_0}\bigr)=
\phi\bigl(\tfrac{c-\mathbf{x}'_t\boldsymbol{\beta}_0}{\sigma_0}\bigr)
)
\end{align}
we have
$E[\mathbf{M}_{11}|\mathbf{x}_t]=E[\mathbf{N}_{11}|\mathbf{x}_t]$.

As to $E[\mathbf{M}_{12}|\mathbf{x}_t]$ and
$E[\mathbf{N}_{12}|\mathbf{x}_t]$,
\begin{align}
&E[\mathbf{M}_{12}|\mathbf{x}_t]\notag\\
&=E\Bigl[(\gamma_0y_t-\mathbf{x}'_t\mathbf{\delta}_0)\mathbf{x}_t
\bigl[\tfrac{1}{\gamma_0}-(\gamma_0y_t-\mathbf{x}'_t\mathbf{\delta}_0)y_t\bigr]\mid
\mathbf{x}_t,y_t>c\Bigr]\texttt{Prob}[y_t>c|\mathbf{x}_t]
+\lambda^2(-v_{0t})(-\mathbf{x}_tc)\notag\\
&\quad\,\cdot \texttt{Prob}[y_t=c|\mathbf{x}_t]\notag
\end{align}
\begin{align}
&=E\Bigl[(\gamma_0y_t-\mathbf{x}'_t\mathbf{\delta}_0)\mathbf{x}_t
\bigl[\tfrac{1}{\gamma_0}-\tfrac{1}{\gamma_0}(\gamma_0y_t-\mathbf{x}'_t\mathbf{\delta}_0)
(\gamma_0y_t-\mathbf{x}'_t\mathbf{\delta}_0)-\tfrac{1}{\gamma_0}(\gamma_0y_t-\mathbf{x}'_t\mathbf{\delta}_0)
\mathbf{x}'_t\mathbf{\delta}_0 \bigr]\mid
\mathbf{x}_t,y_t>c\Bigr]\,\notag\\
&\quad \cdot \texttt{Prob}[y_t>c|\mathbf{x}_t]+\lambda^2(-v_{0t})(-\mathbf{x}_tc)\,\texttt{Prob}[y_t=c|\mathbf{x}_t]\notag\\
&=E\Bigl[(\gamma_0y_t-\mathbf{x}'_t\mathbf{\delta}_0)\tfrac{1}{\gamma_0}
-\tfrac{1}{\gamma_0}(\gamma_0y_t-\mathbf{x}'_t\mathbf{\delta}_0)^3
-\tfrac{1}{\gamma_0}(\gamma_0y_t-\mathbf{x}'_t\mathbf{\delta}_0)^2
\mathbf{x}'_t\mathbf{\delta}_0 \mid
\mathbf{x}_t,y_t>c\Bigr]\mathbf{x}_t\,\texttt{Prob}[y_t>c|\mathbf{x}_t]\notag\\
&\quad+\lambda^2(-v_{0t})(-\mathbf{x}_tc)\,\texttt{Prob}[y_t=c|\mathbf{x}_t]\notag\\
&=\Big\{\tfrac{1}{\gamma_0}\lambda(v_{0t})-\tfrac{1}{\gamma_0}[v^2_{0t}\lambda(v_{0t})+2\lambda(v_{0t})]
-\tfrac{1}{\gamma_0}[v_{0t}\lambda(v_{0t})+1](\gamma_0c-v_{0t})\Big\}\mathbf{x}_t\,\texttt{Prob}[y_t>c|\mathbf{x}_t]\notag\\
&\quad+\lambda^2(-v_{0t})(-\mathbf{x}_tc)\,\texttt{Prob}[y_t=c|\mathbf{x}_t]
\qquad \qquad\quad(\textrm{ by (3.1), (3.3), (3.2) and
}\mathbf{x}'_t\mathbf{\delta}_0=\gamma_0c-v_{0t}\,\,)\notag\\
&=\Big\{-\tfrac{1}{\gamma_0}\lambda(v_{0t})-[v_{0t}\lambda(v_{0t})+1]c+\tfrac{1}{\gamma_0}v_{0t}\Big\}
\mathbf{x}_t\Bigl[1-\Phi\bigl(\tfrac{c-\mathbf{x}'_t\boldsymbol{\beta}_0}{\sigma_0}\bigr)\Bigr]\notag\\
&\quad+\lambda^2(-v_{0t})(-\mathbf{x}_tc)\,\Phi\bigl(\tfrac{c-\mathbf{x}'_t\boldsymbol{\beta}_0}{\sigma_0}\bigr)\notag\\
&=\Big\{-\tfrac{1}{\gamma_0}\lambda(v_{0t})-c+\tfrac{1}{\gamma_0}v_{0t}\Big\}
\mathbf{x}_t\Bigl[1-\Phi\bigl(\tfrac{c-\mathbf{x}'_t\boldsymbol{\beta}_0}{\sigma_0}\bigr)\Bigr]
-v_{0t}\lambda(v_{0t})c\mathbf{x}_t\Bigl[1-\Phi\bigl(\tfrac{c-\mathbf{x}'_t\boldsymbol{\beta}_0}{\sigma_0}\bigr)\Bigr]\notag\\
&\quad+\lambda^2(-v_{0t})(-\mathbf{x}_tc)\,\Phi\bigl(\tfrac{c-\mathbf{x}'_t\boldsymbol{\beta}_0}{\sigma_0}\bigr)\notag\\
&=\Big\{-\tfrac{1}{\gamma_0}\lambda(v_{0t})-c+\tfrac{1}{\gamma_0}v_{0t}\Big\}
\mathbf{x}_t\Bigl[1-\Phi\bigl(\tfrac{c-\mathbf{x}'_t\boldsymbol{\beta}_0}{\sigma_0}\bigr)\Bigr]
-v_{0t}\lambda(-v_{0t})c\mathbf{x}_t\Phi\bigl(\tfrac{c-\mathbf{x}'_t\boldsymbol{\beta}_0}{\sigma_0}\bigr)\notag\\
&\quad+\lambda^2(-v_{0t})(-\mathbf{x}_tc)\,\Phi\bigl(\tfrac{c-\mathbf{x}'_t\boldsymbol{\beta}_0}{\sigma_0}\bigr)
\qquad\qquad\quad \bigl(\textrm{ by
}\lambda(v_{0t})\Bigl[1-\Phi\bigl(\tfrac{c-\mathbf{x}'_t\boldsymbol{\beta}_0}{\sigma_0}\bigr)\Bigr]
=\lambda(-v_{0t})\Phi\bigl(\tfrac{c-\mathbf{x}'_t\boldsymbol{\beta}_0}{\sigma_0}\bigr)\,\, \bigr)\notag\\
&=\Big\{-\tfrac{1}{\gamma_0}\lambda(v_{0t})-c+\tfrac{1}{\gamma_0}v_{0t}\Big\}
\mathbf{x}_t\Bigl[1-\Phi\bigl(\tfrac{c-\mathbf{x}'_t\boldsymbol{\beta}_0}{\sigma_0}\bigr)\Bigr]\notag\\
&\quad+\lambda(-v_{0t})[\lambda(-v_{0t})+v_{0t}](-\mathbf{x}_tc)\,\Phi\bigl(\tfrac{c-\mathbf{x}'_t\boldsymbol{\beta}_0}{\sigma_0}\bigr),
\end{align}
and
\begin{align}
&E[\mathbf{N}_{12}|\mathbf{x}_t]\notag\\
&=E[-y_t\mathbf{x}_t|\mathbf{x}_t,y_t>c]\,\texttt{Prob}[y_t>c|\mathbf{x}_t]+
\lambda(-v_{0t})[\lambda(-v_{0t})+v_{0t}](-\mathbf{x}_tc)\,\texttt{Prob}[y_t=c|\mathbf{x}_t]\notag\\
&=E\bigl[-\tfrac{1}{\gamma_0}(\gamma_0y_t-\mathbf{x}'_t\boldsymbol{\delta}_0
+\mathbf{x}'_t\boldsymbol{\delta}_0)\mathbf{x}_t|\mathbf{x}_t,y_t>c\bigr]\,\texttt{Prob}[y_t>c|\mathbf{x}_t]\notag\\
&\quad+\lambda(-v_{0t})[\lambda(-v_{0t})+v_{0t}](-\mathbf{x}_tc)\,\texttt{Prob}[y_t=c|\mathbf{x}_t]\notag\\
&=E\bigl[-\tfrac{1}{\gamma_0}(\gamma_0y_t-\mathbf{x}'_t\boldsymbol{\delta}_0)-\tfrac{1}{\gamma_0}(\gamma_0c-v_{0t})
|\mathbf{x}_t,y_t>c\bigr]\mathbf{x}_t\,\texttt{Prob}[y_t>c|\mathbf{x}_t]\notag\\
&\quad+\lambda(-v_{0t})[\lambda(-v_{0t})+v_{0t}](-\mathbf{x}_tc)\,\texttt{Prob}[y_t=c|\mathbf{x}_t]\notag\\
&=\Big\{-\tfrac{1}{\gamma_0}\lambda(v_{0t})-c+\tfrac{1}{\gamma_0}v_{0t}\Big\}
\mathbf{x}_t\Bigl[1-\Phi\bigl(\tfrac{c-\mathbf{x}'_t\boldsymbol{\beta}_0}{\sigma_0}\bigr)\Bigr]\notag\\
&\quad+\lambda(-v_{0t})[\lambda(-v_{0t})+v_{0t}](-\mathbf{x}_tc)\,\Phi\bigl(\tfrac{c-\mathbf{x}'_t\boldsymbol{\beta}_0}{\sigma_0}\bigr).\notag\\
&\qquad (\textrm{ by (3.1) })
\end{align}
Therefore
$E[\mathbf{M}_{12}|\mathbf{x}_t]=E[\mathbf{N}_{12}|\mathbf{x}_t]$.
Since $\mathbf{M}_{21}=\mathbf{M}'_{12}$ and
$\mathbf{N}_{21}=\mathbf{N}'_{12}$,
$E[\mathbf{M}_{21}|\mathbf{x}_t]=E[\mathbf{N}_{21}|\mathbf{x}_t]$
also holds.

Finally for $E[\mathbf{M}_{22}|\mathbf{x}_t]$ and
$E[\mathbf{N}_{22}|\mathbf{x}_t]$,
\begin{align}
&E[\mathbf{M}_{22}|\mathbf{x}_t]\notag\\
&=E\Bigl[\bigl[\tfrac{1}{\gamma_0}-(\gamma_0y_t-\mathbf{x}'_t\boldsymbol{\delta}_0)y_t\bigr]^2\mid
\mathbf{x}_t,y_t>c\Bigr]\,\texttt{Prob}[y_t>c|\mathbf{x}_t]+\lambda^2(-v_{0t})c^2\,\texttt{Prob}[y_t=c|\mathbf{x}_t]\notag\\
&=E\Bigl[\bigl[\tfrac{1}{\gamma_0}-\tfrac{1}{\gamma_0}(\gamma_0y_t-\mathbf{x}'_t\boldsymbol{\delta}_0)
(\gamma_0y_t-\mathbf{x}'_t\boldsymbol{\delta}_0)-\tfrac{1}{\gamma_0}(\gamma_0y_t-\mathbf{x}'_t\boldsymbol{\delta}_0)\mathbf{x}'_t\boldsymbol{\delta}_0\bigr]^2\mid
\mathbf{x}_t,y_t>c\Bigr]\,\texttt{Prob}[y_t>c|\mathbf{x}_t]\notag\\
&\quad+\lambda^2(-v_{0t})c^2\,\texttt{Prob}[y_t=c|\mathbf{x}_t]\notag\\
&=E\Bigl[\tfrac{1}{\gamma^2_0}\bigl[1-(\gamma_0y_t-\mathbf{x}'_t\boldsymbol{\delta}_0)^2\bigr]^2-2\tfrac{1}{\gamma^2_0}
\bigl[1-(\gamma_0y_t-\mathbf{x}'_t\boldsymbol{\delta}_0)^2\bigr]
(\gamma_0y_t-\mathbf{x}'_t\boldsymbol{\delta}_0)\mathbf{x}'_t\boldsymbol{\delta}_0\notag\\
&\quad+\tfrac{1}{\gamma^2_0}(\gamma_0y_t-\mathbf{x}'_t\boldsymbol{\delta}_0)^2(\mathbf{x}'_t\boldsymbol{\delta}_0)^2\mid
\mathbf{x}_t,y_t>c\Bigr]\,\texttt{Prob}[y_t>c|\mathbf{x}_t]
+\lambda^2(-v_{0t})c^2\,\texttt{Prob}[y_t=c|\mathbf{x}_t]\notag\\
&=E\Bigl[\tfrac{1}{\gamma^2_0}\bigl[1-(\gamma_0y_t-\mathbf{x}'_t\boldsymbol{\delta}_0)^2\bigr]^2
\mid \mathbf{x}_t,y_t>c\Bigr]\,\texttt{Prob}[y_t>c|\mathbf{x}_t]\notag\\
&\quad-E\Bigl[2\tfrac{1}{\gamma^2_0}
\bigl[1-(\gamma_0y_t-\mathbf{x}'_t\boldsymbol{\delta}_0)^2\bigr]
(\gamma_0y_t-\mathbf{x}'_t\boldsymbol{\delta}_0)\mathbf{x}'_t\boldsymbol{\delta}_0\mid
\mathbf{x}_t,y_t>c\Bigr]
\,\texttt{Prob}[y_t>c|\mathbf{x}_t]\notag\\
&\quad+E\Bigl[\tfrac{1}{\gamma^2_0}(\gamma_0y_t-\mathbf{x}'_t\boldsymbol{\delta}_0)^2(\mathbf{x}'_t\boldsymbol{\delta}_0)^2
\mid \mathbf{x}_t,y_t>c\Bigr]\,\texttt{Prob}[y_t>c|\mathbf{x}_t]\notag\\
&\quad+\lambda^2(-v_{0t})c^2\,\texttt{Prob}[y_t=c|\mathbf{x}_t],
\end{align}
where
\begin{align}
&E\Bigl[\tfrac{1}{\gamma^2_0}\bigl[1-(\gamma_0y_t-\mathbf{x}'_t\boldsymbol{\delta}_0)^2\bigr]^2
\mid \mathbf{x}_t,y_t>c\Bigr]\notag\\
&=\tfrac{1}{\gamma^2_0}-\tfrac{2}{\gamma^2_0}E\bigl[(\gamma_0y_t-\mathbf{x}'_t\boldsymbol{\delta}_0)^2
\mid
\mathbf{x}_t,y_t>c\bigr]+\tfrac{1}{\gamma^2_0}E\bigl[(\gamma_0y_t-\mathbf{x}'_t\boldsymbol{\delta}_0)^4
\mid \mathbf{x}_t,y_t>c\bigr]\notag\\
&=\tfrac{1}{\gamma^2_0}-\tfrac{2}{\gamma^2_0}[v_{0t}\lambda(v_{0t})+1]
+\tfrac{1}{\gamma^2_0}\bigl[v^3_{0t}\lambda(v_{0t})+3v_{0t}\lambda(v_{0t})+3\bigr]
\qquad \qquad\,\,\,\, (\textrm{ by (3.2) and (3.4) })\notag\\
&=\tfrac{1}{\gamma^2_0}\bigl[v^3_{0t}\lambda(v_{0t})+v_{0t}\lambda(v_{0t})+2\bigr],
\end{align}

\begin{align}
&E\Bigl[2\tfrac{1}{\gamma^2_0}
\bigl[1-(\gamma_0y_t-\mathbf{x}'_t\boldsymbol{\delta}_0)^2\bigr]
(\gamma_0y_t-\mathbf{x}'_t\boldsymbol{\delta}_0)\mathbf{x}'_t\boldsymbol{\delta}_0\mid
\mathbf{x}_t,y_t>c\Bigr]\notag\\
&=\tfrac{2}{\gamma^2_0}E\bigl[(\gamma_0y_t-\mathbf{x}'_t\boldsymbol{\delta}_0)|\mathbf{x}_t,y_t>c\bigr]
(\mathbf{x}'_t\boldsymbol{\delta}_0)-\tfrac{2}{\gamma^2_0}
E\bigl[(\gamma_0y_t-\mathbf{x}'_t\boldsymbol{\delta}_0)^3\mid\mathbf{x}_t,y_t>c\bigr]
(\mathbf{x}'_t\boldsymbol{\delta}_0)\notag\\
&=\tfrac{2}{\gamma^2_0}\lambda(v_{0t})(\mathbf{x}'_t\boldsymbol{\delta}_0)
-\tfrac{2}{\gamma^2_0}\bigl[v^2_{0t}\lambda(v_{0t})+2\lambda(v_{0t})\bigr](\mathbf{x}'_t\boldsymbol{\delta}_0)
\qquad \qquad \qquad \qquad \quad (\textrm{ by (3.1) and (3.3) })\notag\\
&=-\tfrac{2}{\gamma^2_0}\bigl[v^2_{0t}\lambda(v_{0t})+\lambda(v_{0t})\bigr](\gamma_0c-v_{0t}),
\qquad \qquad \qquad \qquad \qquad \qquad (\textrm{ by
}\mathbf{x}'_t\boldsymbol{\delta}_0=\gamma_0c-v_{0t}\,\,)\notag\\
\end{align}
and
\begin{align}
&E\Bigl[\tfrac{1}{\gamma^2_0}(\gamma_0y_t-\mathbf{x}'_t\boldsymbol{\delta}_0)^2(\mathbf{x}'_t\boldsymbol{\delta}_0)^2
\mid \mathbf{x}_t,y_t>c\Bigr]\notag\\
&=\tfrac{1}{\gamma^2_0}[v_{0t}\lambda(v_{0t})+1](\gamma_0c-v_{0t})^2.
\qquad \qquad \qquad \qquad(\textrm{ by (3.2) and
}\mathbf{x}'_t\boldsymbol{\delta}_0=\gamma_0c-v_{0t}\,\, )\notag\\
\end{align}
Thus
\begin{align}
&E[\mathbf{M}_{22}|\mathbf{x}_t]\notag\\
&=\tfrac{1}{\gamma^2_0}\bigl[v^3_{0t}\lambda(v_{0t})+v_{0t}\lambda(v_{0t})+2\bigr]\,\texttt{Prob}[y_t>c|\mathbf{x}_t]\notag\\
&\quad \,+\tfrac{2}{\gamma^2_0}\bigl[v^2_{0t}\lambda(v_{0t})+\lambda(v_{0t})\bigr](\gamma_0c-v_{0t})\,\texttt{Prob}[y_t>c|\mathbf{x}_t]\notag\\
&\quad\,+\tfrac{1}{\gamma^2_0}[v_{0t}\lambda(v_{0t})+1](\gamma_0c-v_{0t})^2\,\texttt{Prob}[y_t>c|\mathbf{x}_t]\notag\\
&\quad\,+\lambda^2(-v_{0t})c^2\,\texttt{Prob}[y_t=c|\mathbf{x}_t]\notag\\
&=\tfrac{1}{\gamma^2_0}\bigl[v^3_{0t}\lambda(v_{0t})+v_{0t}\lambda(v_{0t})+2\bigr]\,\texttt{Prob}[y_t>c|\mathbf{x}_t]\notag\\
&\quad\,+\tfrac{1}{\gamma^2_0}\Bigl\{2v^2_{0t}\lambda(v_{0t})+2\lambda(v_{0t})+[v_{0t}\lambda(v_{0t})+1](\gamma_0c-v_{0t})\Bigr\}
(\gamma_0c-v_{0t})\,\texttt{Prob}[y_t>c|\mathbf{x}_t]\notag\\
&\quad\,+\lambda^2(-v_{0t})c^2\,\texttt{Prob}[y_t=c|\mathbf{x}_t]
\qquad \qquad \qquad (\textrm{ by combining the 2nd and the 3rd
terms}\,\,
)\notag\\
&=\tfrac{1}{\gamma^2_0}\bigl[v^3_{0t}\lambda(v_{0t})+v_{0t}\lambda(v_{0t})+2\bigr]\,\texttt{Prob}[y_t>c|\mathbf{x}_t]\notag\\
&\quad\,+\tfrac{1}{\gamma^2_0}\Bigl\{v^2_{0t}\lambda(v_{0t})+2\lambda(v_{0t})+[v_{0t}\lambda(v_{0t})+1]\gamma_0c-v_{0t}\Bigr\}
(\gamma_0c-v_{0t})\,\texttt{Prob}[y_t>c|\mathbf{x}_t]\notag\\
&\quad\,+\lambda^2(-v_{0t})c^2\,\texttt{Prob}[y_t=c|\mathbf{x}_t]
\notag\\
&=\tfrac{1}{\gamma^2_0}\bigl[v^3_{0t}\lambda(v_{0t})+v_{0t}\lambda(v_{0t})+2\bigr]\,\texttt{Prob}[y_t>c|\mathbf{x}_t]\notag\\
&\quad\,+\tfrac{1}{\gamma^2_0}\Bigl\{2\lambda(v_{0t})\gamma_0c+v_{0t}\lambda(v_{0t})\gamma^2_0c^2+
\gamma^2_0c^2-v^3_{0t}\lambda(v_{0t})-2v_{0t}\lambda(v_{0t})-2v_{0t}\gamma_0c+v^2_{0t}\Bigr\}     \,\notag\\
&\quad\,\cdot
\texttt{Prob}[y_t>c|\mathbf{x}_t]+\lambda^2(-v_{0t})c^2\,\texttt{Prob}[y_t=c|\mathbf{x}_t]
\notag\\
&=\tfrac{1}{\gamma^2_0}\Bigl\{2\lambda(v_{0t})\gamma_0c+v_{0t}\lambda(v_{0t})\gamma^2_0c^2+
\gamma^2_0c^2-v_{0t}\lambda(v_{0t})-2v_{0t}\gamma_0c+v^2_{0t}+2\Bigr\}
\Bigl[1-\Phi\bigl(\tfrac{c-\mathbf{x}'_t\boldsymbol{\beta}_0}{\sigma_0}\bigr)\Bigr]\notag\\
&\quad\,+\lambda^2(-v_{0t})c^2\Phi\bigl(\tfrac{c-\mathbf{x}'_t\boldsymbol{\beta}_0}{\sigma_0}\bigr)
\notag\\
&=\tfrac{1}{\gamma^2_0}\Bigl\{2\lambda(v_{0t})\gamma_0c+\gamma^2_0c^2-v_{0t}\lambda(v_{0t})-2v_{0t}\gamma_0c+v^2_{0t}+2\Bigr\}
\Bigl[1-\Phi\bigl(\tfrac{c-\mathbf{x}'_t\boldsymbol{\beta}_0}{\sigma_0}\bigr)\Bigr]\notag\\
&\quad\,+v_{0t}\lambda(v_{0t})c^2\Bigl[1-\Phi\bigl(\tfrac{c-\mathbf{x}'_t\boldsymbol{\beta}_0}{\sigma_0}\bigr)\Bigr]
+\lambda^2(-v_{0t})c^2\Phi\bigl(\tfrac{c-\mathbf{x}'_t\boldsymbol{\beta}_0}{\sigma_0}\bigr)
\notag\\
&=\tfrac{1}{\gamma^2_0}\Bigl\{2\lambda(v_{0t})\gamma_0c+\gamma^2_0c^2-v_{0t}\lambda(v_{0t})-2v_{0t}\gamma_0c+v^2_{0t}+2\Bigr\}
\Bigl[1-\Phi\bigl(\tfrac{c-\mathbf{x}'_t\boldsymbol{\beta}_0}{\sigma_0}\bigr)\Bigr]\notag\\
&\quad\,+\lambda(-v_{0t})[\lambda(-v_{0t})+v_{0t}]c^2\Phi\bigl(\tfrac{c-\mathbf{x}'_t\boldsymbol{\beta}_0}{\sigma_0}\bigr).\\
&\quad\,\, \Bigl(\textrm{ by
}\lambda(v_{0t})\Bigl[1-\Phi\bigl(\tfrac{c-\mathbf{x}'_t\boldsymbol{\beta}_0}{\sigma_0}\bigr)\Bigr]=
\lambda(-v_{0t})\Phi\bigl(\tfrac{c-\mathbf{x}'_t\boldsymbol{\beta}_0}{\sigma_0}\bigr)\,\,
\Bigr)\notag
\end{align}
On the other hand,
\begin{align}
&E[\mathbf{N}_{22}|\mathbf{x}_t]\notag\\
&=E\Bigr[\tfrac{1}{\gamma^2_0}+y^2_t|\mathbf{x}_t,
y_t>c\Bigr]\,\texttt{Prob}[y_t>c|\mathbf{x}_t]+\lambda(-v_{0t})[\lambda(-v_{0t})+v_{0t}]c^2\,\texttt{Prob}[y_t=c|\mathbf{x}_t]\notag\\
&=E\Bigr[\tfrac{1}{\gamma^2_0}+\tfrac{1}{\gamma^2_0}(\gamma_0y_t-\mathbf{x}_t\boldsymbol{\delta}_0
+\mathbf{x}_t\boldsymbol{\delta}_0)^2\mid\mathbf{x}_t,
y_t>c\Bigr]\,\texttt{Prob}[y_t>c|\mathbf{x}_t]+\lambda(-v_{0t})[\lambda(-v_{0t})+v_{0t}]c^2\,\notag\\
&\quad\,\cdot \texttt{Prob}[y_t=c|\mathbf{x}_t]\notag
\end{align}
\begin{align}
&=E\Bigr[\tfrac{1}{\gamma^2_0}+\tfrac{1}{\gamma^2_0}(\gamma_0y_t-\mathbf{x}_t\boldsymbol{\delta}_0
)^2+\tfrac{2}{\gamma^2_0}(\gamma_0y_t-\mathbf{x}_t\boldsymbol{\delta}_0
)(\mathbf{x}_t\boldsymbol{\delta}_0)+\tfrac{1}{\gamma^2_0}(\mathbf{x}_t\boldsymbol{\delta}_0)^2
\mid\mathbf{x}_t,
y_t>c\Bigr]\,\texttt{Prob}[y_t>c|\mathbf{x}_t]\notag\\
&\quad\,+\lambda(-v_{0t})[\lambda(-v_{0t})+v_{0t}]c^2\,\texttt{Prob}[y_t=c|\mathbf{x}_t]\notag\\
&=\tfrac{1}{\gamma^2_0}\Bigl\{1+[v_{0t}\lambda(v_{0t})+1]
+2\lambda(v_{0t})(\mathbf{x}_t\boldsymbol{\delta}_0)
+(\mathbf{x}_t\boldsymbol{\delta}_0)^2\Bigr\}
\,\texttt{Prob}[y_t>c|\mathbf{x}_t]
\notag\\
&\quad\,+\lambda(-v_{0t})[\lambda(-v_{0t})+v_{0t}]c^2\,\texttt{Prob}[y_t=c|\mathbf{x}_t]
\qquad \qquad \qquad (\textrm{ by (3.2) and (3.1)}\,\, )\notag\\
&=\tfrac{1}{\gamma^2_0}\Bigl\{1+[v_{0t}\lambda(v_{0t})+1]
+2\lambda(v_{0t})(\gamma_0c-v_{0t}) +(\gamma_0c-v_{0t})^2\Bigr\}
\Bigl[1-\Phi\bigl(\tfrac{c-\mathbf{x}'_t\boldsymbol{\beta}_0}{\sigma_0}\bigr)\Bigr]\notag\\
&\quad\,+\lambda(-v_{0t})[\lambda(-v_{0t})+v_{0t}]c^2\Phi\bigl(\tfrac{c-\mathbf{x}'_t\boldsymbol{\beta}_0}{\sigma_0}\bigr)\notag\\
&=\tfrac{1}{\gamma^2_0}\Bigl\{2
+2\lambda(v_{0t})\gamma_0c-v_{0t}\lambda(v_{0t})
+\gamma^2_0c^2-2\gamma_0cv_{0t}+v^2_{0t}\Bigr\}
\Bigl[1-\Phi\bigl(\tfrac{c-\mathbf{x}'_t\boldsymbol{\beta}_0}{\sigma_0}\bigr)\Bigr]\notag\\
&\quad\,+\lambda(-v_{0t})[\lambda(-v_{0t})+v_{0t}]c^2\Phi\bigl(\tfrac{c-\mathbf{x}'_t\boldsymbol{\beta}_0}{\sigma_0}\bigr).
\end{align}

Comparing (3.27) with (3.26), we can see that
$E[\mathbf{M}_{22}|\mathbf{x}_t]=E[\mathbf{N}_{22}|\mathbf{x}_t]$.

Since
$E[\mathbf{M}_{ij}|\mathbf{x}_t]=E[\mathbf{N}_{ij}|\mathbf{x}_t],\,\,i=1,2,\,j=1,2$,\,\,
$E[\mathbf{s}(\mathbf{w}_t;\boldsymbol{\delta}_0,\gamma_0)\mathbf{s}(\mathbf{w}_t;\boldsymbol{\delta}_0,\gamma_0)'|
\mathbf{x}_t]$\\
$=-E[\mathbf{H}(\mathbf{w}_t;\boldsymbol{\delta}_0,\gamma_0)|\mathbf{x}_t]$.

Therefore, condition $3$ of Proposition 1.3 is satisfied. \qed

Next we claim
\begin{thm}
For the Tobit Model satisfying Assumption 1', if $\mathbf{
E(}\mathbf{x}_t\mathbf{x}'_t\mathbf{)}$ is nonsingular, then
condition $5$ of Proposition 1.3 is satisfied.
\end{thm}
\pf
 Since condition $3$ of Proposition 1.3 is satisfied,
\begin{equation}
-E[\mathbf{H}(\mathbf{w}_t;\boldsymbol{\delta}_0,
\gamma_0)]=E[\mathbf{s}(\mathbf{w}_t;\boldsymbol{\delta}_0,
\gamma_0) \mathbf{s}(\mathbf{w}_t;\boldsymbol{\delta}_0,
\gamma_0)'].
 \end{equation}

Clearly $\mathbf{s}(\mathbf{w}_t;\boldsymbol{\delta}_0, \gamma_0)
\mathbf{s}(\mathbf{w}_t;\boldsymbol{\delta}_0, \gamma_0)'$ is
positive semidefinite. This means that $\displaystyle
E[\mathbf{s}(\mathbf{w}_t;\boldsymbol{\delta}_0, \gamma_0)\\
\cdot\mathbf{s}(\mathbf{w}_t;\boldsymbol{\delta}_0, \gamma_0)']$
is positive semidefinite. Thus
$E[\mathbf{H}(\mathbf{w}_t;\boldsymbol{\delta}_0, \gamma_0)]$ is
negative semidefinite.

Let
$\mathbf{z}\equiv(\mathbf{z}_1,...,\mathbf{z}_K,\mathbf{z}_{K+1})'\in
\R^{K+1}$ be a solution to the equation
\begin{equation}\label{A}
      {\mathbf{z}}'E[\mathbf{H}(\mathbf{w}_t;\boldsymbol{\delta}_0,
\gamma_0)]\mathbf{z}=0.
\end{equation}

We know that if $\mathbf{z}=(0,...,0)$ is the only solution to the
above equation, then\\
$E[\mathbf{H}(\mathbf{w}_t;\boldsymbol{\delta}_0, \gamma_0)]$ is
nonsingular.

Define $\mathbf{\tilde{A}}$ as the matrix
\begin{equation}
\begin{bmatrix}
\mathbf{x}_t\mathbf{x}'_t&-y_t\mathbf{x}_t\\
-y_t\mathbf{x}'_t&\frac{1}{{\gamma_0}^2}+y^2_t
\end{bmatrix}.
\end{equation}

Define $\mathbf{\tilde{B}}$ as the matrix
\begin{equation}
\begin{bmatrix}
\mathbf{x}_t\mathbf{x}'_t&-c\mathbf{x}_t\\
-c\mathbf{x}'_t&c^2
\end{bmatrix}.
\end{equation}

Then in terms of the expression
$\mathbf{H}(\mathbf{w}_t;\boldsymbol{\delta},\gamma)$ in (1.18),
\begin{equation}
\mathbf{H}(\mathbf{w}_t;\boldsymbol{\delta}_0,
\gamma_0)=-(1-D_t)\mathbf{\tilde{A}}-D_t\lambda(-v_{0t})[\lambda(-v_{0t})+v_{0t}]\mathbf{\tilde{B}}.
\end{equation}

It is easy to see that\footnote{see hint of exercise 3 on page 521
of Hayashi (2000).}
\begin{equation}
\mathbf{\tilde{A}}=\begin{bmatrix}
\mathbf{x}_t\\
-y_t
\end{bmatrix}
\begin{bmatrix}
\mathbf{x}'_t& -y_t
\end{bmatrix}
+ \begin{bmatrix}
\mathbf{0}&\mathbf{0}\\
\mathbf{0'}&\frac{1}{{\gamma_0}^2}
\end{bmatrix},
\end{equation}
and
\begin{equation}
\mathbf{\tilde{B}}=\begin{bmatrix}
\mathbf{x}_t\\
-c
\end{bmatrix}
\begin{bmatrix}
\mathbf{x}'_t& -c
\end{bmatrix}
\end{equation}

Thus combining these with the fact that $1-D_t\geq 0$ and
$D_t\lambda(-v_{0t})[\lambda(-v_{0t})+v_{0t}]\geq 0$, we can see
that $\mathbf{H}(\mathbf{w}_t;\boldsymbol{\delta}_0, \gamma_0)$ is
negative semidefinite.

Suppose $\mathbf{z}_{K+1}\neq 0$, then
\begin{align}
\mathbf{z}'\mathbf{\tilde{A}}\mathbf{z}
&=\biggl[\sum_{i=1}^{K}\mathbf{x}_{ti}\mathbf{z}_i-y_t\mathbf{z}_{K+1}\biggr]^2+\mathbf{z}^2_{K+1}\frac{1}{\gamma^2_0}\notag\\
&\geq \mathbf{z}^2_{K+1}\frac{1}{\gamma^2_0}\notag\\
&>0.
\end{align}

In terms of the expression $D_t$ in (1.15) and the expression
$y^{*}_t$ in (1.10), we have
\begin{align}
 \texttt{Prob} (1-D_t{=}1)&= \texttt{Prob}(y^{*}_t>c) \notag \\
                     &=\texttt{Prob}\bigl(\epsilon_t>c-\tfrac{1}{\gamma_0}\mathbf{x}'_t\boldsymbol{\delta}_0\bigr).
\end{align}

Since $\epsilon_t|\mathbf{x}_t\sim N(0,\sigma^2_0)$ (by (1.11)),
$\texttt{Prob}\bigl(\epsilon_t>c-\tfrac{1}{\gamma_0}\mathbf{x}'_t\boldsymbol{\delta}_0\bigr)>0$.
Thus $\texttt{Prob} (1-D_t{=}1)>0$. Combining this with (3.35), we
have
\begin{equation}
E\bigl[(1-D_t)\mathbf{z}'\mathbf{\tilde{A}}\mathbf{z}\bigr]>0,
\end{equation}
when $\mathbf{z}_{K+1}\neq 0$.

Since $\mathbf{\tilde{B}}$ is positive semidefinite and
$D_t\lambda(-v_{0t})[\lambda(-v_{0t})+v_{0t}]\geq 0$,
\begin{equation}
E\bigl[D_t\lambda(-v_{0t})[\lambda(-v_{0t})+v_{0t}]\mathbf{z}'\mathbf{\tilde{B}}\mathbf{z}\bigr]\geq
0.
\end{equation}

Thus, when $\mathbf{z}_{K+1}\neq 0$,
\begin{align}
{\mathbf{z}}'E[\mathbf{H}(\mathbf{w}_t;\boldsymbol{\delta}_0,
\gamma_0)]{\mathbf{z}}&=
E[{\mathbf{z}}'\mathbf{H}(\mathbf{w}_t;\boldsymbol{\delta}_0,
\gamma_0){\mathbf{z}}]\notag\\
&=E\bigl[(1-D_t)\mathbf{z}'\mathbf{\tilde{A}}\mathbf{z}\bigr]
+E\bigl[D_t\lambda(-v_{0t})[\lambda(-v_{0t})+v_{0t}]\mathbf{z}'\mathbf{\tilde{B}}\mathbf{z}\bigr]\notag\\
&>0,
\end{align}
But this contradicts the assumption
${\mathbf{z}}'E[\mathbf{H}(\mathbf{w}_t;\boldsymbol{\delta}_0,
\gamma_0)]{\mathbf{z}}=0$ (see (3.29)). Therefore
$\mathbf{z}_{K+1}$ must be $0$. This means that
$\mathbf{z}=(\mathbf{z}_1,...,\mathbf{z}_K,0)$, thereafter,
\begin{align}
\mathbf{z}E[\mathbf{H}(\mathbf{w}_t;\boldsymbol{\delta}_0,
\gamma_0)]\mathbf{z}'
=&(\mathbf{z}_1,...,\mathbf{z}_K,0)E[\mathbf{H}(\mathbf{w}_t;\boldsymbol{\delta}_0,
\gamma_0)](\mathbf{z}_1,...,\mathbf{z}_K,0)'\notag\\
=&-(\mathbf{z}_1,...,\mathbf{z}_K)
E\Bigl[\bigl(1-D_t+D_t\lambda(-v_{0t})[\lambda(-v_{0t})+v_{0t}]\bigr)\mathbf{x_tx'_t}\Bigr](\mathbf{z}_1,...,\mathbf{z}_K)'\notag\\
=&-E\Bigl[(\mathbf{z}_1,...,\mathbf{z}_K)
\bigl(1-D_t+D_t\lambda(-v_{0t})[\lambda(-v_{0t})+v_{0t}]\bigr)\mathbf{x_tx'_t}(\mathbf{z}_1,...,\mathbf{z}_K)'\Bigr],
\end{align}

It is clear that for any positive constant
$\overline{\overline{v}} \in \R$, there exists a positive constant
$\overline{\overline{C}}$ which depends on
$\overline{\overline{v}} \in \R$, such that on the set
$\{|v_{0t}|\leq {\overline{\overline{v}}}\}$,
$1-D_t+D_t\lambda(-v_{0t})[\lambda(v_{-0t})-v_{0t}]\geq
\overline{\overline{C}}$. This implies that
\begin{align}
&\,E\Bigl[(\mathbf{z}_1,...,\mathbf{z}_K)
\bigl(1-\lambda(v_{0t})[\lambda(v_{0t})-v_{0t}]\bigr)\mathbf{x_tx'_t}(\mathbf{z}_1,...,\mathbf{z}_K)'\Bigr]\notag\\
\geq&\,E\Bigl[(\mathbf{z}_1,...,\mathbf{z}_K) 1_{\{|v_{0t}|\leq
{\overline{\overline{v}}}\}}\bigl(1-D_t+D_t\lambda(-v_{0t})[\lambda(v_{-0t})-v_{0t}] \bigr)\mathbf{x_tx'_t}(\mathbf{z}_1,...,\mathbf{z}_K)'\Bigr]\notag\\
\geq&\,\overline{\overline{C}}E\Bigl[(\mathbf{z}_1,...,\mathbf{z}_K)
1_{\{|v_{0t}|\leq
{\overline{\overline{v}}}\}}\mathbf{x_tx'_t}(\mathbf{z}_1,...,\mathbf{z}_K)'\Bigr].
\end{align}

The remaining argument follows the same line as we did in the
truncation regression model.

By (3.29), (3.40) and (3.41),
\begin{equation}
{\mathbf{z}}'E[\mathbf{H}(\mathbf{w}_t;\boldsymbol{\delta}_0,
\gamma_0)]\mathbf{z}=0\,\, \Longrightarrow\,\,
(\mathbf{z}_1,...,\mathbf{z}_K)E[ 1_{\{|v_{0t}|\leq
{\overline{v}}\}}\mathbf{x_tx'_t}](\mathbf{z}_1,...,\mathbf{z}_K)'=0.
\end{equation}

It is easy to see that if $E[\mathbf{x_tx'_t}]$ is nonsingular,
then for large enough $\overline{\overline{v}}$,
$E[1_{\{|v_{0t}|\leq
{\overline{\overline{v}}}\}}\mathbf{x_tx'_t}]$ is also
nonsingular. Thus when $E[\mathbf{x_tx'_t}]$ is nonsingular, for
large enough $\overline{\overline{v}}$,
$\{\mathbf{z}_i=0,\,i=1,...,K\}$ is the only solution satisfying
\begin{equation}
(\mathbf{z}_1,...,\mathbf{z}_K)E[ 1_{\{|v_{0t}|\leq
{\overline{\overline{v}}}\}}\mathbf{x_tx'_t}](\mathbf{z}_1,...,\mathbf{z}_K)'=0.
\end{equation}

This means that
\begin{equation}
{\mathbf{z}}'E[\mathbf{H}(\mathbf{w}_t;\boldsymbol{\delta}_0,
\gamma_0)]\mathbf{z}=0\,\, \Longrightarrow\,\,
\mathbf{z}=(0,...,0,0)'.
\end{equation}
Therefore when $E[\mathbf{x_tx'_t}]$ is nonsingular,
$E[\mathbf{H}(\mathbf{w}_t;\boldsymbol{\delta}_0, \gamma_0)]$ is
nonsingular, i.e. condition $5$ of Proposition 1.3 holds. \qed

Theorem 3.2 and Theorem 3.3 imply that the conditional ML
estimator $(\boldsymbol{\hat{\delta}},\hat{\gamma})$ is asymptotic
normal with
$\texttt{Avar}(\boldsymbol{\hat{\delta}},\hat{\gamma})$ given by
the following:
\begin{equation}
\texttt{Avar}(\boldsymbol{\hat{\delta}},\hat{\gamma})=-\{E[\mathbf{H}(\mathbf{w}_t;\boldsymbol{\delta}_0,\gamma_0)]\}^{-1}=\{E[\mathbf{s}(\mathbf{w}_t;\boldsymbol{\delta}_0,\gamma_0)
\mathbf{s}(\mathbf{w}_t;\boldsymbol{\delta}_0,\gamma_0)']\}^{-1}.
\end{equation}

To recover original parameters and obtain the asymptotic variance
of $(\boldsymbol{\hat{\beta}},\hat{\sigma}^2)$, the delta method
can be applied. (see page 520 and page 521 of Hayashi (2000).)

\vspace{.5in}
\begin{singlespace}

\end{singlespace}

\end{doublespace}

\end{document}